\documentclass[a4paper,twoside,12pt]{article}

\usepackage{mathrsfs}
\usepackage{stmaryrd}
\usepackage{amssymb}
\usepackage{amsmath}
\usepackage{exscale}
\usepackage{makeidx}
\usepackage[latin1]{inputenc}
\usepackage{a4}
\usepackage{fancyhdr}
\usepackage[compact,nobottomtitles,rigidchapters]{titlesec}
\usepackage{theorem}
\usepackage{titletoc}

\newcommand\itemi{$(i) \; $}
\newcommand\itemii{$(ii) \; $}
\newcommand\itemiii{$(iii) \; $}
\newcommand\itemiv{$(iv) \; $}

\newcommand\linspan{\text{span}}
\newcommand\closedspan{\overline{\text{span}}}

\newcommand\Sp{\text{Sp}}
\newcommand\Spec{\text{Spec}}

\newcommand\onto{\twoheadrightarrow}
\newcommand\into{\hookrightarrow}

\newcommand\ilim{\varinjlim}
\newcommand\plim{\underset{\longleftarrow}{\lim}\:}

\newcommand{\definition}{\mathrel{\mathop:}=}

\newtheorem{defi}{Definition}
\newtheorem{lemma}{Lemma}
\newtheorem{proposition}{Proposition}
\newtheorem{theorem}{Theorem}

\newtheorem{remark}{Remark}
\newtheorem{observation}{Observation}

\newenvironment{proof}{\paragraph{Proof:}}{\vspace{0.2cm} \hfill $ \boxempty $ \vspace{0.2cm}}

\titleformat{\part}[block]{\Huge\bf}{\thepart. }{0pt}{}
\titlespacing*{\section}{0in}{*4}{*2}
\titlespacing*{\subsection}{0in}{*4}{*2}

\fancypagestyle{plain}{ \fancyhf{} }

\def\romanenum{\renewcommand{\labelenumi}{\textup{(}\roman{enumi}\textup{)}}}

\romanenum

\titlecontents{section}
[0pt]
{}
{\bfseries\contentsmargin{0pt}{\thecontentslabel\enspace}}
{\bfseries\contentsmargin{0pt}}
{\titlerule*[3mm]{.}\contentspage}

\begin{document}

\pagestyle{fancy}

\title{\bfseries The Regular C*-algebra of an Integral Domain}

\author{Joachim Cuntz and Xin Li}

\date{}

\maketitle

\begin{abstract}
To each integral domain $R$ with finite quotients we associate a purely infinite simple C*-algebra in a very natural way.
Its stabilization can be identified with the crossed product of the algebra of continuous functions on the ``finite adele space''
corresponding to $R$ by the action of the $ax+b$-group over the quotient field $Q(R)$.
We study the relationship to generalized Bost-Connes systems and deduce for them a description as universal C*-algebras with the help of our construction.
\end{abstract}

\setcounter{tocdepth}{1}
\tableofcontents

\thanks{\scriptsize{\bf  2000 Math.\ Subject Classification}. Primary: 58B34, 46L05; Secondary: 11R04, 11R56.}

\thanks{\scriptsize{Research supported by the Deutsche Forschungsgemeinschaft (SFB 478).}}

\section{Introduction}

In \cite{Cun1}, the first named author had introduced
C*-algebras $\mathcal{ Q }_\mathbb{Z}$ and $\mathcal{ Q }_\mathbb{N}$ 
associated with the ring of integers $\mathbb{Z}$ or also with the semiring
$\mathbb{N}$, respectively, and which can be obtained from the
natural actions of $\mathbb{Z}$ and $\mathbb{N}$, by multiplication and
addition, on the Hilbert spaces $\ell^2(\mathbb{Z})$ and
$\ell^2(\mathbb{N})$.

This was originally motivated by the well-known construction
by Bost and Connes \cite{BoCo} who had introduced a
C*-dynamical system $(\mathcal{ C } _\mathbb{Q}, \sigma_t)$ and studied
its KMS-states.

The Bost-Connes algebra $\mathcal{ C }_\mathbb{Q}$ is naturally
embedded into $\mathcal{ Q }_\mathbb{N}$. The difference between the two
algebras lies in the fact that $\mathcal{ Q }_\mathbb{N}$ contains,
besides the operators induced by multiplication in $\mathbb{N}$ also
those corresponding to addition.

A main result in \cite{Cun1} was the proof that the algebras
$\mathcal{ Q }_\mathbb{Z}$ and $\mathcal{ Q }_\mathbb{N}$ are purely
infinite simple and, after stabilization, can also be described as
crossed products of the algebra of functions on the finite
adele space $\mathbb{A}_f$ over $\mathbb{Q}$ by the $ax+b$-groups over
$\mathbb{Q}$ or $\mathbb{Q}^+$. This leads in particular to a simple
presentation, by generators and relations, of the C*-algebras
generated by the ``left regular representations'' of $\mathbb{Z}$
and $\mathbb{N}$.

In the present paper we extend the construction of
\cite{Cun1} to an arbitrary commutative ring $R$ without zero
divisors (an integral domain) subject to a finiteness
condition which is typically satisfied by the integral domains
considered in number theory (rings of integers in algebraic number 
fields or polynomial rings over finite fields). We denote the associated
C*-algebra by $\mathfrak{ A } [R]$. We generalize the result from
\cite{Cun1} by showing that $\mathfrak{ A } [R]$ and its
stabilization $\mathfrak{ A } (R)$ are purely infinite simple and
that $\mathfrak{ A } (R)$ can be represented as a crossed
product of the algebra of functions on the ``finite adele space''
corresponding to the profinite completion $\hat{R}$ of $R$,
by the action of the $ax+b$-group over the quotient field
$Q(R)$. At the same time we streamline and improve the
arguments given in \cite{Cun1} in the case $R=\mathbb Z$.

We also show that the higher dimensional analogues of the
Bost-Connes system studied in \cite{CMR} for imaginary
quadratic number fields and in \cite{LLN} for arbitrary
number fields, embed into the C*-algebra $\mathfrak{ A } [R]$ if
the number field allows at most one real place and the class number is 
one. We use this to
deduce a description of the algebras considered in
\cite{CMR}, \cite{LLN} (for number fields with at most one
real place and class number one) in terms of generators and relations.

This description can be used to construct all extremal 
$KMS_\beta$-states of the dynamical system
considered in \cite{LLN} in a very natural way (in complete analogy to 
the original case of $\mathbb{Q}$ treated in \cite{BoCo}).

\section{Universal C*-algebras}

Throughout this article, $ R $ will denote an integral domain with the following properties:

1. the set of units $ R^* $ in $ R $ does not equal $ R^\times \definition R \backslash \left\{ 0 \right\} $ (so we exclude fields)

2. for each $ m \in R^\times $ the ideal $ (m) $ generated by $ m $ in $ R $ is of finite index in $ R $.

We will always think of $ R $ as a subring of its quotient field $ Q(R) $.

Now, let us introduce our C*-algebras $\mathfrak{A}[R]$ in a universal way in terms of generators and relations. Later on, we will see more concrete models for $\mathfrak{A}[R]$. 

\begin{defi}
Let $ \mathfrak{A} [R] $ be the universal C*-algebra generated by isometries $ \left\{ s_m \text{: } m \in R^\times \right\} $ and unitaries $ \left\{ u^n \text{: } n \in R \right\} $ with the relations

\begin{itemize}
  \item[\itemi] $ s_k s_m = s_{km} $
  \item[\itemii] $ u^l u^n = u^{l+n} $
  \item[\itemiii] $ s_m u^n = u^{mn} s_m $
  \item[\itemiv] $ \sum_{ n+(m) \in R / (m) } u^n e_m u^{-n} = 1 $
\end{itemize}

for all $ k, m \in R^\times $, $ l, n \in R $, where $e_m = s_m s_m^*$ is the final projection corresponding to $ s_m $.

\end{defi}

The sum is taken over all cosets $ n+(m) $ in $ R / (m) $ and $ u^n e_m u^{-n} $ is independent of the choice of
$ n $. This follows from \itemi , \itemii and \itemiii (once they are valid).

$ \mathfrak{A} [R] $ exists as the generators must have norm $1$. To show that this universal C*-algebra is not trivial, it suffices to give a non-trivial explicit representation of these generators and relations on a Hilbert space. For this purpose, we consider the ``left regular representation'' on the Hilbert space $ \ell^2 (R) $ (actually, as we are dealing with commutative rings, there is no need to distinguish between ``left'' and ``right'') given by the operators 
\begin{eqnarray}
&& S_m (\xi_r) \definition \xi_{mr} \nonumber \\
&& U^n (\xi_r) \definition \xi_{n+r}. \nonumber
\end{eqnarray}
One immediately checks the relations, \itemiii corresponding to distributivity and \itemiv reflecting
the fact that $ U^n S_m S_m^* U^{-n} $ is the projection onto 
\\
$ \closedspan \left( \left\{ \xi_r \text{: } r \in n+(m) \right\} \right) $ 
and that $ R $ is the disjoint union of the cosets 
\\
$ \left\{ n+(m) \text{: } n \in R \right\} $. 

Therefore, the universal property provides a non-trivial representation via $ s_m \mapsto S_m $, $ u^n \mapsto U^n $.

In analogy to the case of groups, one can think of 
\begin{equation} 
\mathfrak{A}_r [R] \definition C^* \left( \left\{ S_m \text{: } m \in R^\times \right\} \cup \left\{ U^n \text{: } n \in R \right\} \right) \subset \mathcal{L} ( \ell^2 (R) ) \nonumber
\end{equation}
as the reduced (or regular) C*-algebra associated to $R$.

Moreover, denote the $ax+b$-semigroup 
$ \left\{ 
\left(
\begin{smallmatrix}a&b\\0&1\end{smallmatrix}
\right)
\text{ : } a \in R^\times \text{ , } b \in R 
\right\} $
by $ P_R $. 
We have a natural representation of $ P_R $ given by
$
\left(
\begin{smallmatrix}a&b\\0&1\end{smallmatrix}
\right) 
\longmapsto u^b s_a
$.

\section{The Inner Structure}

In order to see that $ \mathfrak{A} [R] $ is simple and purely infinite, we proceed similarly as in \cite{Cun1}. This
means we construct a faithful conditional expectation out of certain group actions and describe this
expectation with the help of appropriate projections (actually, this idea already appears in \cite{Cun2}).

\subsection{Preparations}

We begin with some immediate consequences of the characteristic relations defining $ \mathfrak{A} [R] $. First of
all, the projections $ u^n e_m u^{-n} $, $ u^l e_m u^{-l} $ are orthogonal if $ n+(m) \neq l+(m) $ because
of \itemiv. Denote by $ P $ the set of all these projections, $ P = \left\{ u^n e_m u^{-n} \text{ : } m \in R^\times \text{ , } n \in R \right\} $.
We have the following
\begin{lemma}
\label{immcon}
The formula 
\begin{equation}
e_m = \sum_{ n+(k) \in R / (k) } u^{mn} e_{km} u^{-mn} \nonumber
\end{equation}
is valid for all $ k \text{, } m \in R^\times $.
Furthermore, the projections in $ P $ commute and $ \linspan(P) $ is multiplicatively closed.
\end{lemma}

\begin{proof}
This follows by 
\begin{eqnarray}
& & e_m = s_m 1 s_m^* \nonumber \\
& & = s_m ( \sum_{ n+(k) \in R/(k) } u^n e_k u^{-n} ) s_m^* \nonumber \\
& & = \sum_{ n+(k) \in R/(k) } u^{mn} e_{km} u^{-mn}. \nonumber
\end{eqnarray}
Given two projections $ u^n e_m u^{-n} $, $ u^l e_k u^{-l} $, we can use the formula above to write both projections as
sums of conjugates of $ e_{km} $. Hence it follows that they commute and that their product is in $ \linspan(P) $. 
\end{proof}

As $ \linspan(P) $ is obviously a subspace closed under involution, we get that $ C^*(P) = \closedspan (P) $
is a commutative C*-subalgebra of $ \mathfrak{A} [R] $. We denote it by $ \mathfrak{D}[R] $ and investigate its structure later
on.

Now we present the ``standard form'' of elements in the canonical dense subalgebra of $ \mathfrak{A} [R] $.
\begin{lemma}
\label{densesub}
Set $ S \definition \left\{ s_m^* u^n f u^{-n'} s_{m'} \text{ : } m \text{, } m' \in R^\times \text{ ; } n, n' \in R \text{ ; } f \in P \right\} $. 
\\
Then $ \linspan(S) $ is the smallest *-algebra in $ \mathfrak{A} [R] $ containing the generators 
\\
$ \left\{ s_m \text { : } m \in R^\times \right\} \cup \left\{ u^n \text{ : } n \in R \right\} $.
\end{lemma}
\begin{proof}
Since $ S $ contains the generators and is a subset of the smallest *-algebra containing them, we just have to prove that $ \linspan(S) $ is closed under multiplication
(as it obviously is an involutive subspace). This follows from the following calculation:
\begin{eqnarray}
& & s_m^* u^n f u^{-n'} s_{m'} \cdot s_k^* u^l e u^{-l'} s_{k'} \nonumber \\
&=& s_m^* u^n f u^{-n'} s_{m'} s_{m'}^* s_k^* s_{m'} u^l e u^{-l'} s_{k'} \nonumber \\
&=& s_m^* u^{n-n'} \underbrace{ u^{n'} f u^{-n'} }_{\tilde{f}} e_{m'} s_k^* s_{m'} \underbrace{ u^l e u^{-l} }_{\tilde{e}} u^{l-l'} s_{k'} \nonumber \\
&=& s_m^* u^{n-n'} s_k^* s_k \tilde{f} s_k^* s_k e_{m'} s_k^* s_{m'} \tilde{e} s_{m'}^* s_{m'} u^{l-
l'} s_{k'} \nonumber \\
&=& s_{km}^* u^{kn-kn'} \underbrace{ s_k \tilde{f} s_k^* }_{\in P} \underbrace{ s_k e_{m'} s_k^* }_{\in P} \underbrace{ s_{m'} \tilde{e} s_{m'}^* }_{\in P} u^{lm'-l'm'} s_{k'm'}. \nonumber
\end{eqnarray}
As $ \linspan (P) $ is closed under multiplication, we conclude that the same holds for $ \linspan (S) $.
\end{proof}

\subsection{A Faithful Conditional Expectation}

\begin{proposition}
\label{Ex_faiconexp}
There is a faithful conditional expectation 
\begin{equation}
\Theta \text{: } \mathfrak{A} [R]  \longrightarrow \mathfrak{D}[R] \nonumber
\end{equation}
characterized by 
\begin{equation} 
\Theta ( s_m^* u^n f u^{-n}s_{m'}) = \delta_{m, m'} \delta_{n, n'} s_m^* u^n f u^{-n} s_m \nonumber
\end{equation}
for all $ m $, $ m' \in R^\times $; $ n $, $ n' \in R $; $ f \in P $.
\end{proposition}
\begin{proof}
$ \Theta $ will be constructed as the composition of two faithful conditional expectations
\begin{eqnarray}
&& \Theta_s \text{: } \mathfrak{A} [R] \longrightarrow C^* \left( \left\{ e_m \text{: } m \in R^\times \right\} \cup \left\{ u^n \text{: } n \in R \right\} \right) \nonumber \\
&& \Theta_u \text{: } \Theta_s ( \mathfrak{A} [R] ) \longrightarrow \mathfrak{D}[R] \nonumber
\end{eqnarray}

both arising from group actions on $ \mathfrak{A} [R] $ or $ \Theta_s (  \mathfrak{A} [R]  ) $ respectively.

1. Construction of $ \Theta_s $:

Consider the Pontrjagin dual group $ \hat{G} $ of the discrete multiplicative group $ G \definition (
Q(R)^\times, \cdot) $ in the quotient field of $R$. To each character $ \phi $ in $\hat{G}$ we assign the
automorphism $\alpha_\phi \in Aut( \mathfrak{A} [R] ) $ given by $ \alpha_\phi ( s_m ) = \phi (m) s_m $, $ \alpha_\phi
(u^n) = u^n $ for all $ m \in R^\times $, $ n \in R $. The existence of $\alpha_\phi $ is guaranteed by the
universal property of $ \mathfrak{A} [R] $. In this way, we get a group-homomorphism 
\begin{eqnarray}
  \hat{G} & & \longrightarrow Aut( \mathfrak{A} [R] ) \nonumber \\
  \phi & & \longmapsto \alpha_\phi \nonumber
\end{eqnarray}
which is continuous for the point-norm topology.

It is known that $ \Theta_s $ defined by 
\begin{equation}
\Theta_s (x) = \int_{ \hat{G} } \alpha_\phi (x) d \mu (\phi) \nonumber
\end{equation}
is a faithful conditional expectation from $ \mathfrak{A} [R] $ onto the fixed-point algebra $ \mathfrak{A} [R] ^{ \hat{G} }$, 
where $\mu$ is the normalized Haar measure on the compact group $\hat{G}$ (see \cite{Bla}, II.6.10.4 (v) ). 

It will be useful to determine $ \mathfrak{A} [R] ^{ \hat{G} } $ more precisely. In order to do so let us calculate
\begin{eqnarray}
& & \Theta_s ( s_m^* u^n f u^{-n'} s_{m'} ) \nonumber \\
&=& \int_ { \hat{G} } \alpha_\phi ( s_m^* u^n f u^{-n'} s_{m'} ) d \mu ( \phi ) \nonumber \\
&=& \left( \int_ { \hat{G} } \phi ( m^{-1} m' ) d \mu ( \phi ) \right)  s_m^* u^n f u^{-n'} s_{m'} \nonumber \\
&=& \delta_{m, m'}  s_m^* u^n f u^{-n'} s_{m'}. \nonumber
\end{eqnarray}
Therefore we have 
\\
$ \mathfrak{A} [R] ^{ \hat{G} } = \Theta_s ( \mathfrak{A} [R] ) = \closedspan ( \left\{ s_m^* u^n f u^{-n'}
s_m \text{: } m \in R^\times \text{; } n \text{, } n' \in R \text{; } f \in P \right\} $
\\
as $ \mathfrak{A} [R] = \closedspan ( \left\{ s_m^* u^n f u^{-n'} s_{m'} \text{: } m \text{, } m' \in R^\times \text{; } n \text{, } n' \in R \text{; } f \in P \right\} $ by Lemma \ref{densesub}. But we can even do better claiming 
\\
$ \mathfrak{A} [R] ^{ \hat{G} } = \closedspan ( \left\{ u^n e_m u^{-n'} \text{: } m \in R^\times \text{; } n \text{, } n' \in R \text{; } \right\} ) $,
\\
because we have
\begin{eqnarray}
& & s_m^* u^n e_k u^{-n'} s_m \nonumber \\
&=& s_m^* e_m u^n e_k u^{-n} u^{n-n'} s_m \nonumber \\
&=& s_m^* \sum_{ l + (k) \in R /(k) } u^{lm} e_{km} u^{-lm} u^n \sum_{ i + (m) \in R /(m) } u^{ik}
e_{km} u^{-ik} u^{-n} u^{n-n'} s_m \nonumber \\
&=& s_m^* \sum_a u^{am} e_{km} u^{-am} u^{n-n'} s_m \nonumber \\
&=& \sum_a u^{a} e_{k} u^{-a} s_m^* u^{n-n'} s_m. \nonumber
\end{eqnarray}
where the sums are taken over appropriate indices $ a $ (this being justified by Lemma \ref{immcon}).

Additionally, 
\begin{eqnarray}
& & s_m^* u^{n-n'} s_m = s_m^* e_m u^{n-n'} e_m s_m \nonumber \\
&=& 
  \begin{cases}
    0 \text{ if } n-n' \notin (m) \\
    u^{ m^{-1} ( n-n' ) } \text{ if } n-n' \in (m)
  \end{cases} \nonumber
\end{eqnarray}
so that each $ s_m^* u^n f u^{-n'} s_m $ lies in $ \closedspan ( \left\{ u^n e_m u^{-n'} \text{: } m \in R^\times \text{; } n \text{, } n' \in R \right\} )$. This implies that $ \mathfrak{A} [R] ^{\hat{G}} = C^* ( \left\{ e_m \text{: } m \in R^\times \right\} \cup \left\{ u^n \text{: } n \in R \right\})$.

2. Construction of $\Theta_u$:

Defining $H \definition ( R , + )$ , we have for each $ \chi \in \hat{H} $ an automorphism 
\\
$ \beta_\chi \in Aut ( \Theta_s ( \mathfrak{A} [R] ) ) $ with the properties $ \beta_\chi (e_m) = e_m $ and $ \beta_\chi (u^n) = \chi
(n) u^n$. To see existence of $ \beta_\chi $, we fix $ m \in R^\times $ and consider $C^* ( \left\{ e_m \right\} \cup \left\{ u^n \text{: } n \in R \right\} )$.
\begin{lemma}
\label{lem3}
This algebra is the universal C*-algebra generated by unitaries $ \left\{ v^n \text{: } n \in R \right\} $ and one projection $ f_m $ such that
\begin{eqnarray}
v^n v^{n'} = v^{n+n'} \nonumber \\
\sum_{ n+(m) \in R / (m) } v^n f_m v^{-n} = 1, \nonumber
\end{eqnarray}
the latter relation implicitly including $ v^{lm} f_m = f_m v^{lm} $ for all $l \in R$.
\end{lemma}
\begin{proof}
The universal C*-algebra corresponding to these generators and relations above can be faithfully
represented on a (necessarily infinite-dimensional) Hilbert space. Then it turns out that this
algebra is isomorphic to 
\\
$ M_p ( C^* ( \left\{ v^n \text{: } n \in R \right\} ) ) $ with $ p \definition \# R / (m) $. The isomorphism is provided by the p pairwise orthogonal projections $ v^n f_m v^{-n} $ each being
equivalent to $1$ (where $ \left\{ n + (m) \right\} = R / (m) $). Now the same argument shows 
\\
$ C^* ( \left\{ e_m \right\} \cup \left\{ u^n \text{: } n \in R \right\} ) \cong M_p ( C^* ( \left\{ u^n \text{: } n \in R \right\} ) )$. Thus it remains to show that 
\begin{eqnarray}
C^* ( \left\{ v^n \text{: } n \in R \right\} ) & & \longrightarrow C^* ( \left\{ u^n \text{: } n \in R \right\} ) \nonumber \\
v^n & & \longmapsto u^n \nonumber
\end{eqnarray}
is an isomorphism. This follows by the following observations:

For each n, $ \Sp (u^n) $ is maximal, meaning that it is $\mathbb{T}$ if $ \text{char} (R) = 0 $ and $ \left\{ \zeta \in \mathbb{T} \text{: } \zeta^p = 1 \right\} $ if $ \text{char} (R) = p $ (in this case we have $ (u^n)^p = 1 $ for all $ n \in R $). This follows from the ``left regular representation'' of $ \mathfrak{A} [R] $ discussed above. Therefore, $ \Sp (v^n) = \Sp (u^n) $ for all $ n \in R$.
\\
Given $n_1, ..., n_i \in R$, we have $C^* ( \left\{ v^{n_1} \text{, ... , } v^{n_i} \right\} ) \cong C^* ( \left\{ u^{n_1} \text{, ... , } u^{n_i} \right\} ) $. To see this, we can assume that the $n_1$, ..., $n_i$ are linearly independent over the prime ring of $R$, so that we get 
$ \Spec ( C^* ( \left\{ u^{n_1} \text{, ... , } u^{n_i} \right\} ) $ 
\\
$\cong \Sp( u^{n_1} ) \times ... \times \Sp( u^{n_i} )$
$= \Sp( v^{n_1} ) \times ... \times \Sp( v^{n_i} ) \cong \Spec( C^*( \left\{ v^{n_1} \text{, ... , } v^{n_i} \right\} ) $ which is all we have to show. Now the claim follows by taking the inductive limit of the isomorphisms obtained via the identification of these spectra, and we again get an isomorphism sending $ v^n $ to $u^n $.
\end{proof}

The Lemma shows that we have an automorphism 
\\
$ \beta_{\chi, m} : C^* ( \left\{ e_m \right\} \cup \left\{ u^n \text{: } n \in R \right\} ) 
\longrightarrow C^* ( \left\{ e_m \right\} \cup \left\{ u^n \text{: } n \in R \right\} )$ with
\\ 
$ \beta_{\chi, m} ( e_m ) = e_m $, $ \beta_{\chi, m} ( u^n ) = \chi ( n ) u^n$. 

$\beta_\chi$ can be constructed as the inductive limit of these $\beta_{ \chi , m }$ because 
\\
$ C^* ( \left\{e_m \text{: } m \in R^\times \right\} \cup \left\{ u^n \text{: } n \in R \right\} ) = \ilim C^* ( \left\{ e_m \right\} \cup \left\{ u^n \text{: } n \in R \right\} ) $
with the inclusions 
$ C^* ( \left\{ e_m \right\} \cup \left\{ u^n \text{: } n \in R \right\} ) \into C^* ( \left\{ e_{km} \right\} \cup \left\{ u^n \text{: } n \in R \right\} ) $ 
justified by Lemma \ref{immcon} and as 
$ \beta_{\chi , km} \vert C^* ( \left\{ e_m \right\} \cup \left\{ u^n \text{: } n \in R \right\} ) = \beta_{\chi , m} $. 

Clearly, $\hat{H}$ acts on $\Theta_s ( \mathfrak{A} [R] )$ via $\chi \longmapsto \beta_\chi$ which is again continuous
for the point-norm topology. So we can proceed just as before defining
\begin{equation}
\Theta_u (y) = \int_{ \hat{H} } \beta_\chi ( y) d \mu ( \chi ), \nonumber
\end{equation}
and an analogous calculation shows $ \Theta_u ( u^n e_m u^{-n'} ) = \delta_{n, n'} u^n e_m u^{-n} $. 
Hence it follows that $ ( \Theta_s ( \mathfrak{A} [R] ) )^{ \hat{H} } = \left( \mathfrak{A} [R] ^{ \hat{G} } \right)^{ \hat{H} } = \mathfrak{D}[R] $.

As mentioned at the beginning, we set $ \Theta \definition \Theta_u \circ \Theta_s $ which obviously
yields a faithful conditional expectation with the property

$ \Theta ( s_m^* u^n f u^{-n'} s_{m'} ) = \Theta_u ( \delta_{m, m'} s_m^* u^n f u^{-n'} s_m ) =
\delta_{m , m'} \delta_{n, n'} s_m^* u^n f u^{-n} s_m $.
\end{proof}

In the following we want to give an alternative description of $ \Theta $ with the help of sufficiently
small projections. Let $y$ be in $\linspan (S)$, which means
\begin{equation}
  y = \sum_{ m , m' , n , n' , f } a_{( m , m' , n , n' , f )} s_m^* u^n f u^{-n'} s_{m'}. \nonumber
\end{equation}
In this sum, there are only finitely many projections lying in $P$ which appear with non-trivial coefficients. 
Write them as sums of mutually orthogonal projections $ u^{n_1} e_M u^{-n_1} $, ... , $u^{n_N} e_M u^{-n_N}$.
\begin{proposition}
\label{altproj}
There are $N$ pairwise orthogonal projections $f_i$ in $P$ such that
\begin{itemize}
  \item[I.] $\Phi$ defined by 
  \begin{eqnarray}
  C^* ( \left\{ u^{n_1} e_M u^{-n_1} \text{, ... , } u^{n_N} e_M u^{-n_N} \right\} ) 
  & \longrightarrow & C^* ( \left\{ f_1 \text{, ... , } f_N \right\} ) \nonumber \\
  z & \longmapsto & \sum_{i=1}^N f_i z f_i \nonumber
  \end{eqnarray}
  is an isomorphism.
  \item[II.] $ \Phi ( \Theta (y) ) = \sum_{i=1}^N f_i y f_i $
\end{itemize}
\end{proposition}
\begin{proof}
We will find appropriate $\nu_i$ and $\mu$ so that $f_i \definition u^{\nu_i} e_{\mu} u^{-\nu_i}$
satisfies I. and II.

As a first step, the conditions
\begin{eqnarray}
&& \nu_i + (M) = n_i + (M) \text{ for all } 1 \leq i \leq N \nonumber \\
&& \mu \in (M) \nonumber
\end{eqnarray}
enforce mutual orthogonality and imply I. as we have for $\lambda = M^{-1} \mu$ (in $R$ by the second
condition)
\begin{eqnarray}
& & f_i u^{n_j} e_M u^{-n_j} f_i \nonumber \\
&=& f_i \sum_{l + (\lambda) \in R / (\lambda) } u^{n_j + lM} e_\mu u^{-n_j - lM} f_i \nonumber \\
&=& \delta_{i,j} f_i \nonumber
\end{eqnarray}
because 
\begin{eqnarray}
& & f_i u^{n_j + lM} e_\mu u^{-n_j - lM} \neq 0 \text{ for some } l \in R \nonumber \\
&\Leftrightarrow& \nu_i + (\mu) = n_j + lM + (\mu) \text{ for some } l \in R \nonumber \\
&\Leftrightarrow& \nu_i + (M) = n_j + (M) \nonumber \\
& & (\mu \in (M)) \nonumber \\
&\Leftrightarrow& i = j \nonumber \\
& & (\nu_i + (M) = n_i + (M) \neq n_j + (M) \text{ for all } i \neq j). \nonumber
\end{eqnarray}
Therefore, $\Phi$ maps $u^{n_i} e_M u^{-n_i}$ to $f_i$ and is thus an isomorphism.

To find sufficient conditions on $\nu_i$ and $\mu$ for II., let us consider those summands in $y$ with
$a_{(m , m' , n , n' , f)} \neq 0$ and $\delta_{m , m'} \delta_{n , n'} = 0$. Call the corresponding
indices $(m , m' , n , n', f)$ critical, there are only finitely many of them. As $\Theta$ maps such
summands to $0$, we have to ensure that $f_i s_m^* u^n f u^{-n'} s_{m'} f_i = 0$ for those critical
indices, and as $ \Theta $ acts identically on summands with $\delta_{m, m'} \delta_{n, n'} = 1$, this will be sufficient for II.

We have 
\begin{eqnarray}
& & f_i s_m^* u^n f u^{-n'} s_{m'} f_i \nonumber \\
&=& s_m^* u^n \left( u^{-n} s_m f_i s_m^* u^n \right) f \left( u^{-n'} s_{m'} f_i s_{m'}^* u^{n'} \right) u^{-n'} s_{m'} \nonumber \\
&=& s_m^* u^n ( u^{m \nu_i - n} e_{m \mu} u^{-m \nu_i + n} u^{m' \nu_i - n'} e_{m' \mu} u^{-m' \nu_i +
n'} ) f u^{-n'} s_{m'} \nonumber 
\end{eqnarray}
and the term in brackets can be described as
\begin{eqnarray}
& & u^{m \nu_i - n} e_{m \mu} u^{-m \nu_i + n} u^{m' \nu_i - n'} e_{m' \mu} u^{-m' \nu_i + n'} \nonumber \\
&=& \sum_{ a + (m') \in R / (m') } u^{-n+m \nu_i + a m \mu} e_{m m' \mu} u^{n-m \nu_i - a m \mu} \nonumber \\
&\cdot& \sum_{ b + (m) \in R / (m) } u^{-n'+m' \nu_i + b m' \mu} e_{m m' \mu} u^{n-m' \nu_i - b m' \mu}. \nonumber
\end{eqnarray}
Now we see that the projections in these two sums are pairwise orthogonal if 
\begin{eqnarray}
& & -n + m \nu_i + a m \mu + (m m' \mu) \neq -n' + m' \nu_i + b m' \mu + (m m' \mu) \nonumber \\
&\Leftrightarrow& n-n' + \nu_i (m'-m) + (b m' - a m) \mu \notin ( m m' \mu) \text{ for all } a \text{, } b \text{ in } R. \nonumber
\end{eqnarray}
This can be enforced by the even stronger condition
\begin{equation}
n-n' + \nu_i (m'-m) \notin (\mu), \nonumber
\end{equation}
which we have to satisfy for each critical index simultaneously.

On the whole, the projections $f_i$ satisfy I. and II. if $\nu_i$ and $\mu$ have the three properties
\begin{itemize}
  \item $ \nu_i + (M) = n_i + (M) $ for all $ 1 \leq i \leq N $
  \item $ \mu \in (M) $
  \item $ n-n' + \nu_i (m'-m) \notin (\mu) $ for all critical indices.
\end{itemize}
One could, for example, choose $\nu_i$ such that $\nu_i + (M) = n_i + (M)$ for all $1 \leq i \leq N$ and
$n-n' + \nu_i (m'-m) \neq 0$ for all critical indices. This can be simultaneously done as there are
infinitely many possibilities for the $\nu_i$ to satisfy the first condition, while the second one only
excludes finitely many (namely $- (m'-m)^{-1} (n-n')$ for all critical indices with $m \neq m'$,
otherwise this condition is automatically valid as $\delta_{m , m'} \delta_{n , n'} = 0$). Then just
take an element $r \in R^\times $ which is not invertible and set 
\begin{equation}
\mu \definition r M \prod [ n-n' + \nu_i (m'-m) ] \in R^\times \nonumber
\end{equation}
where the product is taken over all critical indices. It is immediate that this choice of $\mu$ enforces the second and third condition.
\end{proof}

\subsection{Purely Infinite Simple C*-algebras}

With the help of these ingredients it is now possible to prove the following result:
\begin{theorem}
$ \mathfrak{A} [R] $ is simple and purely infinite, i.e. for all $ 0 \neq x \in \mathfrak{A} [R] $ there are $a$, $b$
$\in \mathfrak{A} [R] $ with $a x b = 1$. 
\end{theorem}
\begin{proof}
Consider first a positive, non-trivial element $x$ in $ \mathfrak{A} [R] $. Recall that we have constructed a faithful conditional expectation $ \Theta $ in Proposition \ref{Ex_faiconexp}. As $\Theta (x) \neq 0$ we can assume
$\lVert \Theta (x) \rVert = 1$. As $ \linspan (S)$ is dense in $ \mathfrak{A} [R] $ (compare Lemma \ref{densesub}) we can find $y \in \linspan (S)_+$ with $\lVert x-y \rVert < \frac{1}{2}$, $\lVert \Theta (y) \rVert = 1$. Proposition \ref{altproj} gives us pairwise orthogonal
projections $f_i$ and $ \Phi $ depending on $y$ such that 
\begin{equation}
\Phi ( \Theta (y) ) = \sum_{i=1}^N f_i y f_i = \sum_{i=1}^N \lambda_i f_i \nonumber
\end{equation}
for some non-negative $\lambda_i$ as we know that $\Phi ( \Theta (y) )$ lies in $C^*( \left\{ f_1 \text{, ... , } f_N \right\} )$ and that $\Theta (y)$ is positive. Since $\Phi$ is isometric, we have $1 = \lVert \Phi ( \Theta (y) ) \rVert$, so that there must be an index $j$ with $\lambda_j = 1$, as $\lVert \sum_{i=1}^N \lambda_i f_i \rVert = \sup_{1 \leq i \leq N} \lambda_i$. Consider the isometry $s \definition u^{ \nu_j } s_\mu$. It has the properties
\begin{eqnarray}
&& s s^* = f_j \text{ and } s^* f_j s = s^* s s^* s = 1 \text{ so that} \nonumber \\
&& s^* y s = s^* f_j s s^* y s s^* f_j s = s^* f_j^2 y f_j^2 s = s^* f_j s = 1. \nonumber
\end{eqnarray}
Therefore, we conclude that
\begin{equation}
\lVert s^* x s - 1 \rVert = \lVert s^* (x-y) s \rVert < \frac{1}{2} \nonumber
\end{equation}
which implies that $s^* x s$ is invertible in $ \mathfrak{A} [R] $.

Set $a \definition ( s^* x s )^{-1} s^*$ and $b \definition s$, this gives $a x b = ( s^* x s )^{-1} s^*
x s = 1$ as claimed.

Given an arbitrary non-trivial element $x$, we get by the same argument as above, used on $x^* x$,
elements $a'$ and $b'$ with $a' x^* x b' = 1$ so that we can set $a \definition a' x^*$ and $b \definition
b'$.
\end{proof}
\begin{remark}
\label{remuniv}
An immediate consequence is the fact that every C*-algebra generated by unitaries and isometries
satisfying the characteristic relations is canonically isomorphic to $ \mathfrak{A} [R] $. 

As a special case of this observation, we get $ \mathfrak{A}_r [R] \cong \mathfrak{A} [R] $.
\end{remark}

\section{Representation as a Crossed Product}

This section is about representing $\mathfrak{A}[R]$ as a crossed product involving some kind of a generalized finite adele ring and the $ax+b$-group $P_{Q(R)}$.

\subsection{Ring-theoretical Constructions}

We start with some ring-theoretical constructions. Set
\begin{equation}
\hat{R} = \plim \left\{ R/(m) \text{; } p_{m,lm} \right\} \nonumber
\end{equation}
where $ p_{m,lm} \text{: } R/(lm) \longrightarrow R/(m)$ is the canonical projection. This is the profinite
completion of $R$.

A concrete description would be 
\begin{equation}
\hat{R} = \left\{ (r_m)_m \in \prod_{m \in R^\times} R/(m) : p_{m,lm} (r_{lm}) = r_m \right\} \nonumber
\end{equation}
with the induced topology of the product $\prod_m R/(m)$, each finite Ring $R/(m)$ carrying the discrete topology. $\hat{R}$ is a compact ring with addition and multiplication defined componentwise. Furthermore, we have
the diagonal embedding 
\begin{eqnarray}
R && \into \hat{R} \nonumber \\
r && \longmapsto (r)_m \nonumber
\end{eqnarray}
and we will identify $R$ with a subring of $\hat{R}$ via this embedding.

Moreover, for $l \in R^\times$ we have the canonical projection $\hat{R} \onto R / (l)$. Its kernel
equals $l \hat{R}$ as those elements are apparently mapped to $0$, while an element $(r_m)_m \in \hat{R}$
mapped to $0$ can be written as $l \cdot ( l^{-1} r_{lm} )_m \in l \hat{R}$. Therefore, we get an
isomorphism $\hat{R} / l \hat{R} \cong R / (l)$.

As a next step, set
\begin{equation}
\label{adelic}
\mathscr{R} \definition \ilim \left\{ \mathscr{R}_m \text{; } \phi_{m, lm} \right\} \nonumber
\end{equation}
where $\mathscr{R}_m = \hat{R}$ for all $m \in R^\times$ and $\phi_{m, lm}$ is multiplication with $l$.
\\
An explicit picture for $\mathscr{R}$ is 
\begin{equation} 
\coprod_{m \in R^\times} \mathscr{R}_m / \sim \nonumber
\end{equation}
where $x_l \sim y_m$ 
$ \Leftrightarrow m x_l = l y_m \text{ for } x_l \in \mathscr{R}_l \text{, } y_m \in \mathscr{R}_m$. Denote by $p$ the canonical projection $\coprod_{m \in R^\times} \mathscr{R}_m \onto \mathscr{R} $ and by $\iota_m$ the embedding
\begin{eqnarray}
\hat{R} && \rightarrow \mathscr{R}_m \rightarrow \mathscr{R} \nonumber \\
x && \longmapsto p(x). \nonumber
\end{eqnarray}
$\mathscr{R}$ is a locally compact ring via 
\begin{eqnarray}
&& \iota_m (x) + \iota_l (y) = \iota_{lm} ( lx + my ), \nonumber \\
&& \iota_m (x) \cdot \iota_l (y) = \iota_{lm} (xy). \nonumber
\end{eqnarray}
Again, we identify $\hat{R}$ with a subring of $\mathscr{R}$ via $\iota_1$.
\\
An immediate observation is the fact that $\iota_m ( \hat{R} )$ is compact and open in $\mathscr{R}$.
Compacity is clear as $\iota_m$ is continuous and $\hat{R}$ is compact. Furthermore, 
\begin{eqnarray}
&& \mathscr{R}_l \cap p^{-1} ( \iota_m ( \hat{R} ) ) \nonumber \\
&=& \left\{ x_l \in \mathscr{R}_l \text{: } x_l \sim y_m \text{ for some } y_m \in \mathscr{R}_m \right\} \nonumber \\
&=& \phi^{-1}_{l, lm} ( l \hat{R} ) \nonumber
\end{eqnarray}
and $l \hat{R}$ is open in $\hat{R}$ because $\hat{R} \backslash l \hat{R} = \bigcup_{  r + (l) \neq (l)} r + l \hat{R}$
using the isomorphism $\hat{R} / l \hat{R} \cong R / (l) $, so that $\hat{R} \backslash l \hat{R}$ is a finite 
union of compact sets, thus closed.

\subsection{Description of the Algebra}

With these preparations, we can establish connections with the C*-algebra $ \mathfrak{A} [R] $.
\begin{observation}
\label{firstob}
We have $\mathfrak{D}[R] \cong C( \hat{R} )$ via $ u^n e_m u^{-n} \longmapsto p_{m \hat{R} + n} $, where $p_{m \hat{R} + n}$ denotes the characteristic function on the compact and open subset $m \hat{R} + n \subset \hat{R}$.
\end{observation}
\begin{proof}
$\mathfrak{D}[R]$ can be described as the inductive limit of 
\\
$D_m = C^* ( \left\{ u^n e_m u^{-n} \text{: } n \in R / (m) \right\} )$ with the inclusions $ D_m \into D_{lm} $. Furthermore, $ \Spec(D_m) \cong R/(m) $ as the projections $u^n e_m u^{-n}$
are mutually orthogonal, and 
\begin{eqnarray}
Spec (D_{lm}) && \longrightarrow Spec (D_m) \nonumber \\
\chi && \longmapsto \chi \vert D_m \nonumber
\end{eqnarray}
corresponds to 
\begin{eqnarray}
p_{lm,m} \text{: } R / (lm) && \longrightarrow R / (m) \nonumber \\
r+(lm) && \longmapsto r+(m) \nonumber
\end{eqnarray}
via this identification. Therefore, we have $ \Spec (D) \cong \plim \left\{ R/(m) \text{; } p_{m,lm} \right\} = \hat{R} $. Thus we get the isomorphism
\begin{eqnarray}
\alpha \text{: } \mathfrak{D}[R] && \longrightarrow C( \hat{R} ) \nonumber \\
u^n e_m u^{-n} && \longmapsto p_{m \hat{R} + n}. \nonumber
\end{eqnarray}
\end{proof}
\begin{defi}
The stabilization of $ \mathfrak{A} [R] $, denoted by $ \mathfrak{A} (R) $, is defined as the inductive limit of the system
$ \left\{  \mathfrak{A} (R) _m \text{; } \varphi_{m,lm} \right\} $ 
\\
where $ \mathfrak{A} (R)_m = \mathfrak{A} [R] $ and $ \varphi_{m,lm} \text{: } \mathfrak{A} [R] \longrightarrow \mathfrak{A} [R] $ is given by $ x \longmapsto s_l x s_l^* $.

Furthermore, we set $\mathfrak{D}(R) = \ilim \left\{\mathfrak{D}(R)_m \text{; } \varphi_{m,lm} \right\} $ with
$ \mathfrak{D}(R)_m = \mathfrak{D}[R] $ and $\varphi_{m,lm}$ just defined as above. $\mathfrak{D}(R)$ can obviously be identified with a C*-
subalgebra of $ \mathfrak{A} (R) $.
\end{defi}
\begin{observation}
\label{secob}
We have $ \mathfrak{D}(R) \cong C_0 ( \mathscr{R} )$.
\end{observation}
\begin{proof}
The maps $\varphi_{m,lm}$, conjugated by $\alpha$ (see Observation \ref{firstob}), give maps
\begin{equation}
\psi_{m,lm} \definition \alpha \circ \varphi_{m,lm} \circ \alpha^{-1} \text{: }
C ( \hat{R} ) \longrightarrow C ( \hat{R} ) \nonumber
\end{equation}
where $\psi_{m,lm} (f) (x) = f ( l^{-1} x ) p_{l \hat{R} } (x)$. This follows from the calculation
\begin{eqnarray}
&& \psi_{m,lm} \circ \alpha ( u^n e_m u^{-n} ) (x) \nonumber \\
&=& \psi_{m,lm} (p_{m \hat{R} + n}) (x) \nonumber \\
&=& p_{m \hat{R} + n} ( l^{-1} x ) p_{l \hat{R}} (x) \nonumber \\
&=& p_{ l m \hat{R} + l n } (x) \nonumber \\
&=& \alpha ( u^{ln} e_{lm} u^{-ln} ) (x) \nonumber \\
&=& \alpha \circ \varphi_{m,lm} ( u^n e_m u^{-n} ) (x). \nonumber
\end{eqnarray}
This yields an isomorphism $\overline{ \alpha } \text{: } \mathfrak{D}(R) \longrightarrow \ilim \left\{ C( \hat{R} ) \text{; }
\psi_{m,lm} \right\}$. Additionally, we consider homomorphisms 
\begin{eqnarray}
\kappa_k \text{: } C( \hat{R} ) && \longrightarrow C_0 ( \mathscr{R} ) \nonumber \\
f && \longmapsto f \circ \iota_k^{-1} \cdot p_{ \iota_k ( \hat{R} ) }. \nonumber
\end{eqnarray}
They satisfy $ \kappa_{lm} \circ \varphi_{m,lm} = \kappa_m $ because
\begin{eqnarray}
&& \kappa_{lm} \circ \varphi_{m,lm} (f) (x) \nonumber \\
&=& \varphi_{m,lm} (f) ( \iota_{lm}^{-1} (x) ) p_{ \iota_{lm} ( \hat{R} ) } (x) \nonumber \\
&=& f ( l^{-1} \iota_{lm}^{-1} (x) ) \cdot p_{ \iota_{lm} ( l \hat{R} ) } (x) \nonumber \\
&=& f ( \iota_m^{-1} (x) ) p_{ \iota_m } (x) \nonumber \\
&=& \kappa_m (f) (x). \nonumber
\end{eqnarray}
Hence these homomorphisms give a homomorphism 
\begin{equation}
\ilim \left\{ C( \hat{R} ) \text{; } \psi_{m,lm} \right\} \longrightarrow C_0 ( \mathscr{R} ) \nonumber
\end{equation}
which is injective as each $\kappa_k$ is injective because of $\kappa_k (f) \circ \iota_k = f$
and surjective which follows from the fact that $\mathscr{R} = \bigcup_{m \in R^\times} \iota_m ( \hat{R} )$
and Stone-Weierstrass.
\end{proof}

Finally, we come to the already mentioned picture of $ \mathfrak{A} (R) $.

\begin{theorem}
\label{theo2}
$ \mathfrak{A} (R) $ is isomorphic to $C_0 ( \mathscr{R} ) \rtimes P_{ Q(R) }$ where the $ax+b$-group acts on $\mathscr{R}$ via 
affine transformations.
\end{theorem}
\begin{proof}
The first step is the observation that we have a canonical isomorphism 
\begin{equation}
p_{ \hat{R} } ( C_0 (\mathscr{R}) \rtimes P_{ Q(R) } ) p_{ \hat{R} } \cong \mathfrak{A} [R] \nonumber
\end{equation}
denoted by $\beta$:

To this end, consider 
$u^n \definition V_{
\left(
\begin{smallmatrix}1&n\\0&1\end{smallmatrix}
\right)
}
p_{ \hat{R} }$ and $s_m \definition 
V_{
\left(
\begin{smallmatrix}m&0\\0&1\end{smallmatrix}
\right)
} p_{\hat{R} } $. One checks that these are unitaries and isometries in $p_{ \hat{R} } ( C_0 (\mathscr{R})
\rtimes P_{ Q(R) } ) p_{ \hat{R} }$ satisfying the characteristic relations of $ \mathfrak{A} [R] $. Furthermore, we have $u^n
s_m s_m^* u^{-n} = p_{m \hat{R} + n}$ so that $C^*( \left\{ u^n s_m s_m^* u^{-n} \text{: } m \text{, } m' \in R^\times \text{; } n \text{, } n' \in R \right\} )$ is a closed C*-subalgebra of $C( \hat{R} )$ seperating points and thus equal to $C(\hat{R} ) = p_{ \hat{R} } C_0 ( \mathscr{R} )$ by Stone-Weierstrass. Hence it follows that 
$p_{ \hat{R} } ( C_0 (\mathscr{R}) \rtimes P_{ Q(R) } ) p_{ \hat{R} }$ is the C*-algebra generated by the
$u^n$ and $s_m$ and thus isomorphic to $ \mathfrak{A} [R] $ by Remark \ref{remuniv}.

Secondly, define 
\begin{equation}
\tilde{ \varphi }_{m,lm} \text{: } p_{ \hat{R} } ( C_0 (\mathscr{R}) \rtimes P_{ Q(R) } ) p_{ \hat{R} } 
\longrightarrow p_{ \hat{R} } ( C_0 (\mathscr{R}) \rtimes P_{ Q(R) } ) p_{ \hat{R} } \nonumber
\end{equation}
to be conjugation by 
$
V_{
\left(
\begin{smallmatrix}l&0\\0&1\end{smallmatrix}
\right)
} 
$. It is clear that we have $ \beta \circ \varphi_{m,lm} \circ \beta^{-1} = \tilde{ \varphi }_{m,lm}$, thus an isomorphism 
\begin{equation}
\overline{ \beta } \text{: } \mathfrak{A} (R) \longrightarrow 
\ilim \left\{ p_{ \hat{R} } ( C_0 (\mathscr{R}) \rtimes P_{ Q(R) } ) p_{ \hat{R} } ; \tilde{ \varphi }_{m,lm} \right\}. 
\nonumber
\end{equation}
Moreover, set 
\begin{eqnarray}
\lambda_k \text{: } p_{ \hat{R} } ( C_0 (\mathscr{R}) \rtimes P_{ Q(R) } ) p_{ \hat{R} } && \longrightarrow 
C_0 (\mathscr{R}) \rtimes P_{ Q(R) } \nonumber \\
z && \longmapsto 
V_{
\left(
\begin{smallmatrix}k&0\\0&1\end{smallmatrix}
\right)
}^* z 
V_{
\left(
\begin{smallmatrix}k&0\\0&1\end{smallmatrix}
\right)
}. \nonumber
\end{eqnarray}
As 
\begin{eqnarray}
&& \lambda_{lm} \circ \tilde{ \varphi }_{m,lm} (z) \nonumber \\
&=& V_{
\left(
\begin{smallmatrix}lm&0\\0&1\end{smallmatrix}
\right)
}^* 
V_{ 
\left(
\begin{smallmatrix}l&0\\0&1\end{smallmatrix}
\right)
} 
z 
V_{
\left(
\begin{smallmatrix}l&0\\0&1\end{smallmatrix}
\right)
}^* 
V_{
\left(
\begin{smallmatrix}lm&0\\0&1\end{smallmatrix}
\right)
} \nonumber \\
&=& 
V_{
\left(
\begin{smallmatrix}m&0\\0&1\end{smallmatrix}
\right)
}^* z 
V_{
\left(
\begin{smallmatrix}m&0\\0&1\end{smallmatrix}
\right)
} = \lambda_m (z), \nonumber
\end{eqnarray}
this gives a homomorphism
\begin{equation}
\lambda \text{: } \ilim \left\{ p_{ \hat{R} } ( C_0 (\mathscr{R}) \rtimes P_{ Q(R) } ) p_{ \hat{R} } \text{; } \tilde{
\varphi }_{m,lm} \right\} \longrightarrow C_0 (\mathscr{R}) \rtimes P_{ Q(R) } \nonumber
\end{equation}
which is injective as this is the case for each $\lambda_m$, and it is surjective as 
\\
$\lambda_m ( p_{ \hat{R} } ) = p_{ \iota_m ( \hat{R} ) }$ is an approximate unit for $C_0 ( \mathscr{R} ) \rtimes P_{ Q(R) }$.
\end{proof} 

\begin{remark}
Combining this result with the preceding remark, we see that 
\begin{equation}
\mathfrak{A}_r [R] \cong \mathfrak{A} [R] \cong p_{ \hat{R} } ( C_0 (\mathscr{R}) \rtimes P_{ Q(R) } ) p_{ \hat{R} }, \nonumber
\end{equation}
which yields a faithful (and very natural) representation of 
\\
$ p_{ \hat{R} } ( C_0 (\mathscr{R}) \rtimes P_{ Q(R) } ) p_{ \hat{R} } $ on $ \ell^2 (R) $.
\end{remark}

\begin{remark}
We call $ \mathfrak{A} (R) $ the stabilization because $ \mathfrak{A} (R) \cong
\mathcal{K} \otimes \mathfrak{A} [R] $. This comes from the observation that $ \mathfrak{A} [R]
$ is isomorphic to $ M_L ( \mathfrak{A} [R] ) $ with regard to the $ L $ pairwise
orthogonal projections $ \left\{ u^n e_l u^{-n} : n \in R \right\} $ where 
$ L = \# R / (l) $. 
And under this identification, conjugation
with $ s_l $ (which is $\varphi_{m,lm}$) corresponds to the inclusion of $ \mathfrak{A} [R] $ into the upper 
left corner of $ M_L ( \mathfrak{A} [R] ) $. 

In other words, using the theory of crossed products by semigroups, we can also say that 
$ \mathfrak{A} [R] \cong C( \hat{R} ) \rtimes P_R $ and that the dynamical system corresponding to 
$ \mathfrak{A} (R) $ is just the associated minimal dilation system (see \cite{Lac}).
\end{remark}

\begin{remark}
Having the classification programme for C*-algebras in mind, one should note that each of the algebras $ \mathfrak{A} [R] $ is
nuclear as $ P_{Q(R)} \cong Q(R) \rtimes Q(R)^\times $ is always amenable because it is solvable.
\end{remark}

\section{Links to Algebraic Number Theory}

The typical examples we have in mind are the rings of integers in an
algebraic number field and polynomial rings with coefficients in a finite field. 
These are exactly the objects of interest in algebraic number theory.

Let $ R = \mathfrak{o} $ be such a ring (the conditions from the beginning are
satisfied). First of all, we have in this case $ \hat{ \mathfrak{o} } \cong \prod \mathfrak{o}_{\nu} $, 
where the product is taken over all the finite places $ \nu $ over $ K = Q(\mathfrak{o}) $. Here, $ \hat{ \mathfrak{o} } $ is to be understood in the sense of the previous section.

If $ \mathfrak{o} $ has positive
characteristic (i.e. $ \mathfrak{o} $ sits in a finite extension of $ \mathbb{F}_p (T) $), we call the place corresponding
to $ T^{-1} $ infinite. $ \mathfrak{o}_{\nu} $ is the maximal compact subring in the
completion of $ \mathfrak{o} $ with regard to $ \nu $. 

Furthermore, we have $ \mathscr{R} \cong \mathbb{A}_{f, K} $ which is the
finite adele ring; with the notation in \cite{Wei}, IV § 1, this is $ K_A ( \left\{
\text{``infinite places''} \right\} ) $. 
This can be seen as follows: 
\begin{eqnarray}
  & & \hat{\mathfrak{o}} = \plim_{m} \left\{ \mathfrak{o} / (m) \right\} \nonumber \\
  & \cong & \plim_{\wp_i, n_i} \left\{ \mathfrak{o} / \wp_{1}^{n_1} ... \wp_{l}^{n_l} \right\} \nonumber \\
  & & \text{(} \mathfrak{o} \text{ is a Dedekind ring with unique factorization of ideals)} \nonumber \\
  & \cong & \plim_{n_\wp} \left\{ \prod_{ \wp \in Spec( \mathfrak{o} ) \backslash \left\{ 0 \right\} } \mathfrak{o} / \wp^{n_\wp} \right\} \nonumber \\
  & & \text{(Chinese remainder theorem)} \nonumber \\
  & \cong & \prod_{ \wp \in Spec( \mathfrak{o} ) \backslash \left\{ 0 \right\} } \plim_n \left\{ \mathfrak{o} / \wp^n \right\}   
  \nonumber \\
  & \cong & \prod_{ \nu \text{ finite} } \plim_n \left\{ \mathfrak{o}_{\nu} / P_{\nu}^n \right\} \nonumber \\
  & & \text{(there is a bijection between non-trivial prime ideals and finite places)} \nonumber \\
  & \cong & \prod_{ \nu \text{ finite} } \mathfrak{o}_{\nu} \nonumber \\
  & & \text{(} \mathfrak{o} \text{ is a Dedekind ring)} \nonumber
\end{eqnarray}
where $ P_{\nu} $ is the subset of $ \mathfrak{o}_{\nu} $ with valuation strictly
smaller than $ 1 $. 

The second identification comes from
\begin{equation}
\mathscr{R} \cong \ilim \hat{\mathfrak{o}} \cong \ilim \left\{ \prod \mathfrak{o}_{\nu} \right\} \cong ( \mathfrak{o}^\times )^{-1} \prod \mathfrak{o}_{\nu} \cong K_A ( \left\{ \text{``infinite places''} \right\}. \nonumber
\end{equation}
The details can be found in \cite{Wei} and \cite{Neu}.

So all in all, we have purely infinite simple C*-algebras $ \mathfrak{A} [ \mathfrak{o} ] \cong C( \hat{ \mathfrak{o} } ) \rtimes P_{ \mathfrak{o} } $ with
stabilization $ \mathfrak{A} (\mathfrak{o}) \cong C_0 ( \mathbb{A}_{f, K} ) \rtimes P_K $.

\section{Relationship to Generalized Bost-Connes Systems}

As mentioned at the beginning, our investigations are partly motivated
by the work of Bost and Connes, who studied a $C^*$-dynamical system for
$\mathbb{Q}$ which had several interesting properties: e.g. it revealed
connections to explicit class field theory over the rational numbers (see \cite{BoCo}). As
a next step, Connes, Marcolli and Ramachandran succeeded in constructing
a $ C^* $-dynamical system for imaginary quadratic number fields and
establishing analogous connections to explicit class field theory for
these (see \cite{CMR}).

In the meantime, there have been several attempts to construct systems
with similar properties for arbitrary number fields (see \cite{CoMa} for an overview).

Most recently, Laca, Larsen and Neshveyev considered $ C^* $-dynamical
systems for arbitrary number fields generalizing the systems mentioned 
above for the case of $ \mathbb{Q} $ and imaginary quadratic fields. 
Moreover, they managed to classify the corresponding KMS-states, 
which was a key ingredient in setting up connections to class field theory. 
Still, these results have not yet led to more insights concerning explicit class field theory.

Our aim in the following section is to embed these generalized Bost-Connes 
algebras into $ \mathfrak{A} $, at least for a certain class of
number fields.
Viewing these generalized Bost-Connes systems as subalgebras of our $
C^* $-algebra $ \mathfrak{A} $, it will be possible the deduce for them a
description as universal $ C^* $-algebras with generators and relations,
as Bost and Connes originally did in the case of $ \mathbb{Q} $.

Before we compare the $ C^* $-algebras constructed by Laca,
Larsen and Neshveyev with our universal $C^*$-algebras, let us very briefly explain their construction, to
set up the notation:

Fix an algebraic number field K and let $ \mathfrak{o} $ be its ring of
integers.

Denote the ring of finite adeles by 
$ \mathbb{A}_f = \underset{\nu \text{ finite} }{\prod^\prime} K_{\nu} $, 
where we take the restricted direct product 
with respect to the ring inclusions $ \mathfrak{o}_{\nu} \subset K_{\nu} $; 
let $ K_{\infty} = \underset{\nu \text{ infinite} }{\prod} K_{\nu} $ be the product
of infinite places; 
then the ring of adeles can be written as 
$ \mathbb{A} = K_{\infty} \times \mathbb{A}_f $.

Furthermore, the group of ideles is 
$ \mathbb{A}^* = K_{\infty}^{\times} \times \underset{\nu \text{ finite} }{\prod^\prime} K_{\nu}^{\times} $, 
this time the restricted product is taken with respect to 
$ \mathfrak{o}_{\nu}^* \subset K_{\nu}^{\times} $. 
Let us write $ K_{\infty , +} $ for the component of the identity in $ K_{\infty} $.

Moreover, take $ \hat{\mathfrak{o}} = \underset{ \nu \text{ finite} }{\prod} \mathfrak{o}_{\nu} $ and 
$ \hat{\mathfrak{o}}^* = \underset{ \nu \text{ finite} }{\prod} \mathfrak{o}_{\nu}^* $.

We will frequently think of subsets of $ \mathbb{A}_f $ as embedded in $
\mathbb{A} $, just by filling in zeros at the infinite places (or
identities in the multiplicative case). Moreover, the algebraic number
field (or subsets in $K$) can always be thought of as subsets of the
adeles (or of the ideles in the multiplicative case) using the diagonal embedding.

Now, for each number field $K$, Laca, Larsen and Neshveyev define a topological space 

\begin{equation}
X = \mathbb{A}^* / \overline{ K^{\times} K_{\infty , +}^{\times} }
\times_{ \hat{ \mathfrak{o} }^* } \mathbb{A}_f \nonumber
\end{equation}

which is a quotient of $ \mathbb{A}^* / \overline{ K^{\times} K_{\infty , +}^{\times} }
\times \mathbb{A}_f $ with respect to the equivalence relation

\begin{eqnarray}
&& ( (x_{\nu}), (y_{\nu}) ) \sim ((x_{\nu}^\prime), (y_{\nu}^\prime) ) \nonumber \\
&\Leftrightarrow& \text{ there exists } (r_{\nu})
\in \hat{\mathfrak{o}}^* \text{ with } ( (r_{\nu}) (x_{\nu}),
(r_{\nu}) ^{-1} (y_{\nu}) ) = ( (x_{\nu}^\prime), (y_{\nu}^\prime) ). \nonumber
\end{eqnarray}

For brevity, let us write $ U $ for $ \overline{ K^{\times} K_{\infty , +}^{\times} } $.
There is a clopen subset 
\\
$ Y = \mathbb{A}^* / U \times_{ \hat{ \mathfrak{o} }^* } \hat{ \mathfrak{o} } $ 
sitting in $ X $.
Furthermore, they consider an action of 
$ \mathbb{A}_f^* / \hat{ \mathfrak{o} }^* $ on $ X $ given by 
$ (z_{ \nu }) ((x_{ \nu }), (y_{ \nu }) ) = ( (z_{ \nu })^{-1} (x_{\nu}) , (z_{ \nu }) (y_{ \nu }) )$.
Finally, their $C^*$-algebra is given by

\begin{equation}
\mathcal{A} = 1_Y \left( C_0 (X) \rtimes \mathbb{A}_f^* / \hat{ \mathfrak{o} }^* \right) 1_Y. \nonumber
\end{equation}

At this point, we should note that - presented in this way - this is a purely adelic-idelic way of describing the system, but that these objects have their natural meaning in number theory via certain abstract
identifications (mostly provided by class field theory), for instance:

$ \mathbb{A}^* / U \cong Gal(K^{ab} / K) $,
where $ Gal(K^{ab} / K) $ is the Galois group of the maximal abelian
field extension of $K$, or

$ \mathbb{A}_f^* / \hat{ \mathfrak{o} }^* \cong J_K $, 
where $ J_K $ is the group of fractional ideals viewed as a discrete group
(see \cite{Wei}, IV § 3).

\subsection{Comparison of the Adelic-Idelic Constructions}

We start the comparison on a purely topological level considering the
adelic-idelic constructions. The first aim will be to establish a relationship
between
$ \mathbb{A}^* / U \times_{ \hat{ \mathfrak{o} }^* } \mathbb{A}_f $ and 
$ \mathbb{A}_f $.

There is a canonical map
\begin{eqnarray}
\psi^* : \mathbb{A}_f \longrightarrow \mathbb{A}^* / U \times_{
\hat{ \mathfrak{o} }^* } \mathbb{A}_f \nonumber \\
( y_{\nu} ) \longmapsto ( (1)^{\bullet}, (y_{\nu}) )^{\bullet} \nonumber
\end{eqnarray}

which we would like to investigate in detail.

From the definitions, we immediately get
\begin{eqnarray}
&& \psi^* ( (y_{\nu}) ) = \psi^* ( ( \tilde{y}_{\nu} ) ) \nonumber \\
& \Leftrightarrow & ( (1)^{\bullet}, (y_{\nu}) ) \sim ( (1)^{\bullet}, ( \tilde{y}_{\nu} ) ) \nonumber \\ 
& \Leftrightarrow & \text{there exists } (z_{\nu}) \in \hat{
\mathfrak{o} }^* \text{ such that } ( (1)^{\bullet}, (y_{\nu}) ) = (
(z_{\nu})^{\bullet} , (z_{\nu})^{-1} ( \tilde{y}_{\nu} ) ) \nonumber \\
& \Leftrightarrow & \text{there exists } (z_{\nu}) \in \hat{
\mathfrak{o} }^* \cap U \text{ with } (y_{\nu}) = (z_{\nu})^{-1} ( \tilde{y}_{\nu} ). \nonumber
\end{eqnarray}

Let us calculate $ \hat{ \mathfrak{o}}^* \cap U $, as this will be needed later on:

\begin{lemma}
\label{calculationintersection}
\begin{equation}
\hat{ \mathfrak{o}}^* \cap U = \overline{ \mathfrak{o}^* \cap \bigcap_{\nu \text{ real} } \nu^{-1} ( \mathbb{R}_{>0} ) } \nonumber
\end{equation}
\end{lemma}

\begin{proof}
The inclusion ``$\supset$'' holds because we have 
\begin{equation} 
\mathfrak{o}^* \subset \hat{ \mathfrak{o} }^* \text{ and } 
\mathfrak{o}^* \cap \bigcap_{\nu \text{ real}} \nu^{-1} ( \mathbb{R}_{>0} ) \subset K^{\times} K_{\infty , +}^{\times}. \nonumber
\end{equation}

The get the other inclusion, observe
\begin{eqnarray}
&& (z_{\nu}) \in \hat{ \mathfrak{o} }^* \cap U \nonumber \\
& \Leftrightarrow & (z_{\nu}) \in \hat{ \mathfrak{o} }^* \text{ and there
exists a sequence } (z_{\nu}^{ (n) }) \text{ in } K^{\times} K_{\infty ,
+}^{\times} \text{ with } \nonumber \\
&& (z_{\nu}) = \lim_{n \rightarrow \infty}
(z_{\nu}^{ (n) }) \text{ in } \mathbb{A}^*. \nonumber
\end{eqnarray}

By the definition of the topology on $ \mathbb{A}^* $, there is a finite
set of places P such that
\begin{eqnarray}
&& (z_{\nu}) \in \prod_{\nu \in P} K_{\nu}^{\times} \times \prod_{\nu \notin P}
\mathfrak{o}_{\nu}^* \nonumber \\
& \Rightarrow & \text{ there is } \tilde{N} \in \mathbb{N} \text{ with
} (z_{\nu}^{ (n) }) \in \prod_{\nu \in P} K_{\nu}^{\times} \times \prod_{\nu
\notin P} \mathfrak{o}_{\nu}^* \text{ for all } n \geq \tilde{N}. \nonumber
\end{eqnarray}

As $ \underset{n \rightarrow \infty}{\lim} (z_{\nu}^{ (n) }) = (z_{\nu}) $, we
conclude that $ \underset{n \rightarrow \infty}{\lim} z_{\nu}^{ (n) } = z_{\nu} $
for all places $ \nu $ in $ K_{\nu}^{\times} $, but as $ z_{\nu}^{ (n) }
\in \mathfrak{o}_{\nu}^* $ for almost all finite places if $ n \geq \tilde{N} $ and
because $ \mathfrak{o}_{\nu}^* $ is open in $ K_{\nu}^{\times} $, there must be a $
N \in \mathbb{N} $ ($ N \geq \tilde{N} $) such that:
\begin{equation}
z_{\nu}^{ (n) } \in \mathfrak{o}_{\nu}^* \text{ for all finite places } \nu \text{
and for all } n \geq N. \nonumber
\end{equation}
But as $ (z_{\nu}^{ (n) }) $ lies in $ K^{\times} K_{\infty ,
+}^{\times} $, there exists for each $n$ $z^{ (n) }$ in $K^{\times}$
such that $ z_{\nu}^{ (n) } = z^{ (n) } $ at every finite place ($K$ is
diagonally embedded). Thus, we have $ z^{ (n) } = z_{\nu}^{ (n) } \in
K^{\times} \cap \mathfrak{o}_{\nu}^* $ for all finite places and $ n \geq N $,
which means that 
$ z^{ (n) } \in K^{\times} \cap \underset{\nu \text{ finite}}{\bigcap} \mathfrak{o}_{\nu}^* 
= \mathfrak{o}^* $ 
for all $ n \geq N $. Hence it follows that 
\begin{equation}
(z_{\nu}^{ (n) }) = (z^{ (n) }) \in \mathfrak{o}^* \subset \hat{ \mathfrak{o} }^* \text{ for all } n \geq N. \nonumber
\end{equation}

Moreover, as $ (z_{\nu}) $ lies in $ \hat{ \mathfrak{o} }^* $, we must
have $ z_{\nu} = 1 $ for all infinite places, so that $ z_{\nu}^{ (n) }
\in \mathbb{R}_{>0} $ for all real places and sufficiently large $n$.
Using the observation above ($ z^{ (n) } = z_{\nu}^{
(n) } $), we conclude that $ z^{ (n) } \in \nu^{-1} ( \mathbb{R}_{>0} )
$ for all real places and $n$ large enough, and hence 

\begin{equation}
(z_{\nu}) \in
\overline{ \mathfrak{o}^* \cap \bigcap_{\nu \text{ real} } \nu^{-1} (
\mathbb{R}_{>0} ) } \nonumber
\end{equation}
as claimed.

\end{proof}

Let us now consider the quotient space $ \mathbb{A}_f / \sim_{\hat{
\mathfrak{o} }^* \cap U } $, where 
\begin{eqnarray}
&& (y_{\nu}) \sim_{\hat{ \mathfrak{o} }^* \cap U } ( \tilde{y}_{\nu} )
\text{ if and only if there exists } (z_{\nu}) \in \hat{ \mathfrak{o} }^*
\cap U \nonumber \\
&& \text{such that } (y_{\nu}) = (z_{\nu}) (\tilde{y}_{\nu}). \nonumber
\end{eqnarray}

Using the universal property of this quotient, we get a continuous and injective map
\begin{equation}
\varphi^* : \mathbb{A}_f / \sim_{ \hat{ \mathfrak{o} }^* \cap U }
\longrightarrow \mathbb{A}^* / U \times_{ \hat{ \mathfrak{o} }^* }
\mathbb{A}_f \nonumber
\end{equation}

with $ \psi^* = \varphi^* \circ p $, where $p$ is the projection map $
\mathbb{A}_f \rightarrow \mathbb{A}_f / \sim_{\hat{ \mathfrak{o} }^*
\cap U } $.

\begin{lemma}
$ \varphi^* $ is closed.
\end{lemma}

\begin{proof}
It suffices to show that $ \psi^* $ is closed, because of the following:
Take $ A \subset \mathbb{A}_f / \sim_{\hat{ \mathfrak{o} }^* \cap U } $
to be an arbitrary closed set. As $p$ is continuous, $p^{-1} (A) $ is
closed in $ \mathbb{A}_f $. Assuming that $ \psi^* $ is closed, we get
that $ \varphi^* (A) = \varphi^* p p^{-1} (A) = \psi^* p^{-1} (A) $ is
closed in $ \mathbb{A}^* / \sim_{ \hat{ \mathfrak{o} }^* \cap U } $.

Now take $ A \subset \mathbb{A}_f $ to be an arbitrary closed set, and let 
\begin{equation}
\pi : \mathbb{A}^* / U \times \mathbb{A}_f \rightarrow
\mathbb{A}^* / U \times_{ \hat{ \mathfrak{o} }^* } \mathbb{A}_f \nonumber
\end{equation}
be the canonical projection. We have to show that $ \psi^* (A) $ is closed,
which is equivalent to closedness of $ \pi^{-1} \psi^* (A) $. We have
\begin{eqnarray}
&& \pi^{-1} \psi^* (A) \nonumber \\
&=& \left\{ ( (a_{\nu})^{\bullet} , (b_{\nu}) ) \in \mathbb{A}^* \times
\mathbb{A}_f : \pi ( (a_{\nu})^{\bullet}, (b_{\nu}) ) \in \psi^* (A) \right\}
\nonumber \\
&=& \left\{ ( (a_{\nu})^{\bullet} , (b_{\nu}) ) \in \mathbb{A}^* \times
\mathbb{A}_f : \exists \text{ } (y_{\nu}) \in A : 
( (a_{\nu})^{\bullet} , (b_{\nu}) ) \sim ( (1)^{\bullet},(y_{\nu}) ) \right\}
\nonumber \\
&=& \left\{ ( (a_{\nu})^{\bullet} , (b_{\nu}) ) : \exists \text{ } (y_{\nu}) \in A \text{, }
(z_{\nu}) \in \hat{ \mathfrak{o} }^* : 
(a_{\nu})^{\bullet} = (z_{\nu})^{\bullet} \wedge (b_{\nu}) = (z_{\nu})^{-1} (y_{\nu}) \right\}
\nonumber \\
&=& \left\{ ( (z_{\nu})^{\bullet} , (z_{\nu})^{-1} (y_{\nu}) ) \in \mathbb{A}^* \times
\mathbb{A}_f : (z_{\nu}) \in \hat{ \mathfrak{o} }^* , (y_{\nu}) \in A \right\}
\nonumber
\end{eqnarray}
Now suppose we have a sequence 
$ ( (z_{\nu}^{ (n) } )^{\bullet}, (z_{\nu}^{ (n) })^{-1} (y_{\nu}^{ (n) } ) ) $ 
in $ \pi^{-1} \varphi^* (A) $ converging
to $ ( (a_{\nu})^{\bullet}, (b_{\nu}) ) \in \mathbb{A}^* / U \times
\mathbb{A}_f $. Then we claim:
$ ( (a_{\nu})^{\bullet}, (b_{\nu}) ) \in \pi^{-1} \psi^* (A) $.

Indeed: $ \hat{ \mathfrak{o} }^* $ is compact, therefore there is a
convergent subsequence $ (z_{\nu}^{ (n_k) }) $ with limit $ (z_{\nu}) \in
\hat{ \mathfrak{o} }^* $. Thus, $ (z_{\nu})^{\bullet} = \underset{k \rightarrow
\infty}{\lim} (z_{\nu}^{(n_k)})^{\bullet} = \underset{n \rightarrow \infty}{\lim} (z_{\nu}^{ (n) })^{\bullet} =
(a_{\nu})^{\bullet} $ and $ (y_{\nu}^{ (n_k) }) = ( z_{\nu}^{ (n_k) } ) (
z_{\nu}^{ (n_k) } )^{-1} (y_{\nu}^{ (n_k) }) \underset{k}{\longrightarrow} (z_{\nu})
(b_{\nu}) $. Hence, $ (y_{\nu}^{ (n_k) }) $ is a sequence in $A$
converging in $ \mathbb{A}_f $, therefore $ (z_{\nu}) (b_{\nu}) =
\lim_{k \rightarrow \infty} (y_{\nu}^{ (n_k) }) =: (y_{\nu}) $ lies in $A$. 
Thus we have $ (b_{\nu}) = (z_{\nu})^{-1} (y_{\nu}) $ and hence
$ ( (a_{\nu})^{\bullet}, (b_{\nu}) ) = ( (z_{\nu})^{\bullet}, (z_{\nu})^{-1}
(y_{\nu}) ) \in \pi^{-1} \psi^* (A) $ which proves our claim and
therefore the Lemma.

\end{proof}

It remains to investigate under which conditions $ \varphi^* $ is
surjective. 
\begin{lemma}
\label{surjektiv?}
$ \varphi^* $ is surjective if $ h_K = 1 $ and there is at most one real
(infinite) place of $K$.
\end{lemma}
\begin{proof}
First of all, $ \varphi^* $ is surjective if and anly if $ \psi^* $ is
so. Now, $ \psi^* $ is surjective if for any $ (a_{\nu})^{\bullet} \in
\mathbb{A}^* / U $, $ (b_{\nu}) \in \mathbb{A}_f $ there are $ (y_{\nu})
\in \mathbb{A}_f $, $ (z_{\nu}) \in \hat{ \mathfrak{o} }^* $ such that $
(a_{\nu})^{\bullet} = (z_{\nu})^{\bullet} $ and $ (b_{\nu}) = (z_{\nu})^{-1}
(y_{\nu}) $ (which is equivalent to $ ( (a_{\nu})^{\bullet}, (b_{\nu}) ) \sim
( (1)^{\bullet}, (y_{\nu}) ) $). As the first condition is the crucial one
(once it holds, the second one can be enforced), $ \psi^* $ is
surjective if (and only if) $ \mathbb{A}^* = \hat{ \mathfrak{o} }^* \cdot
U $.

Assuming that the number of real places is not bigger then one, we have
$ K^{\times} K_{\infty}^{\times} \hat{ \mathfrak{o} }^* = K^{\times}
K_{\infty , +}^{\times} \hat{ \mathfrak{o} }^* $ because given $ (a)
(b_{\nu}) (c_{\nu}) \in K^{\times} K_{\infty}^{\times} \hat{ \mathfrak{o}
}^* $ and we have $ b_{\nu} < 0 $ at the real place (otherwise there is
nothing to prove), then $ (a) (b_{\nu}) (c_{\nu}) = (-a) (-b_{\nu}) (-
c_{\nu}) \in K^{\times} K_{\infty , +}^{\times} \hat{ \mathfrak{o} }^* $.

Furthermore, $ h_K = 1 $ implies $ 1 = \# J_K / P_K = \# I(K) / P(K) =
\# \mathbb{A}^* / K^{\times} K_{\infty}^{\times} \hat{ \mathfrak{o} }^* $
(see \cite{Wei}, V § 3). Hence it follows that we have $ \mathbb{A}^* =
K^{\times} K_{\infty}^{\times} \hat{ \mathfrak{o} }^* = K^{\times}
K_{\infty , +}^{\times} \hat{ \mathfrak{o} }^* $ which implies that $
\psi^* $ is surjective.

\end{proof}

Summarizing our observations to this point, we get
\begin{proposition}
If $h_K = 1$ and there is at most one real place, then the map
\begin{eqnarray}
\varphi^* : \mathbb{A}_f / \sim_{ \hat{ \mathfrak{o} }^* \cap U } &&
\longrightarrow \mathbb{A}^* / U \times_{ \hat{ \mathfrak{o} }^* }
\mathbb{A}_f \nonumber \\
(y_{\nu})^{\bullet} && \longmapsto ( (1)^{\bullet}, (y_{\nu}) )^{\bullet} \nonumber
\end{eqnarray}
is a homoemorphism.
\end{proposition}

\begin{remark}
One should note that Lemma \ref{surjektiv?} is not optimal in the sense that $ \varphi^* $ can be surjective even if $ K $ has more than one real place. The crucial criterion is whether $ \mathfrak{o}^* $ is embedded in $ K_{\infty}^{\times} $ in such a way that every possible combination of signs (in the real places) can be arranged (compare \cite{LaFr}, Proposition 4 for a similar problem).
\end{remark}

\subsection{Comparison of the C*-algebras}

As a next step, let us study the corresponding $C^*$-algebras and the
crossed products in the situation of the last proposition (we assume $h_K = 1$ and that there is at
most one real place):
\begin{proposition}
\label{algebraiso}
$ \varphi^* $ induces *-isomorphisms
\begin{equation}
C_0 (X) \rtimes \mathbb{A}_f^* / \hat{ \mathfrak{o} }^* 
\overset{\simeq}{\longrightarrow}
C_0 ( \mathbb{A}_f / \sim_{ \hat{ \mathfrak{o} }^* \cap U } ) \rtimes
K^{\times} / \mathfrak{o}^* \nonumber
\end{equation}
and 
\begin{equation}
\mathcal{A} 
\overset{\simeq}{\longrightarrow }
1_{ \hat{ \mathfrak{o} } / \sim_{ \hat{ \mathfrak{o} }^* \cap U } }
\left(
C_0 ( \mathbb{A}_f / \sim_{ \hat{ \mathfrak{o} }^* \cap U } ) \rtimes
K^{\times} / \mathfrak{o}^* 
\right)
1_{ \hat{ \mathfrak{o} } / \sim_{ \hat{ \mathfrak{o} }^* \cap U } } \nonumber
\end{equation}
if there are no real places of $K$. 

If there is a real place, we get *-isomorphisms
\begin{equation}
C_0 (X) \rtimes \mathbb{A}_f^* / \hat{ \mathfrak{o} }^* 
\overset{\simeq}{\longrightarrow}
C_0 ( \mathbb{A}_f / \sim_{ \hat{ \mathfrak{o} }^* \cap U } ) \rtimes
K^{\times}_{>0} / \mathfrak{o}^*_{>0} \nonumber
\end{equation} 
and 
\begin{equation}
\mathcal{A} 
\overset{\simeq}{\longrightarrow}
1_{ \hat{ \mathfrak{o} } / \sim_{ \hat{ \mathfrak{o} }^* \cap U } }
\left(
C_0 ( \mathbb{A}_f / \sim_{ \hat{ \mathfrak{o} }^* \cap U } ) \rtimes K^{\times}_{>0} / \mathfrak{o}^*_{>0} 
\right)
1_{ \hat{ \mathfrak{o} } / \sim_{ \hat{ \mathfrak{o} }^* \cap U } }. \nonumber
\end{equation}

\end{proposition}

Here we have fixed a real embedding corresponding to the real place and we think of $K$ as a subset of $\mathbb{R}$ via this embedding.

\begin{proof}
As a first step, $ \varphi^* $ induces a *-isomorphism 
\begin{equation}
C_0 (X) 
\overset{\simeq}{\longrightarrow} 
C_0 ( \mathbb{A}_f / \sim_{ \hat{ \mathfrak{o} }^* \cap U } ) \nonumber
\end{equation} 
(recall $ X = \mathbb{A}^* / U \times_{\hat{ \mathfrak{o} }^*} \mathbb{A}_f $). 

Now, let us assume that $K$ does not have any real places. It remains to
prove that the actions of $ \mathbb{A}_f^* / \hat{ \mathfrak{o} }^* $ and $
K^{\times} / \mathfrak{o}^* $ are compatible. But this follows from the
fact that we have an isomorphism (because of $h_K = 1$)
\begin{eqnarray}
K^{\times} / \mathfrak{o}^* && \longrightarrow \mathbb{A}^*_f / \hat{
\mathfrak{o} }^* \nonumber \\
\lambda^{\bullet} && \longmapsto (\lambda)^{\bullet} \nonumber
\end{eqnarray}
and the following computation:
\begin{eqnarray}
&& (\lambda)^{\bullet} \cdot \varphi^* ( (y_{\nu}) ) \nonumber \\
& = & (\lambda)^{\bullet} \cdot( (1)^{\bullet}, (y_{\nu}) )^{\bullet} \nonumber \\
& = & ( ( \lambda^{-1} )^{\bullet}, (\lambda) (y_{\nu}) )^{\bullet} \nonumber \\
& = & ( (1)^{\bullet}, (\lambda) (y_{\nu}) )^{\bullet} \nonumber \\
& = & \varphi^* ( \lambda^{\bullet} \cdot (y_{\nu})^{\bullet} ). \nonumber
\end{eqnarray}

This gives us the first isomorphism, which we denote by $ \varphi $. 

To get the second one, we have to show
$ \varphi ( 1_Y ) = 1_{ \hat{ \mathfrak{o} } / \sim_{ \hat{ \mathfrak{o} }^* \cap
U } } $, which is equivalent to $ \varphi^* ( \hat{ \mathfrak{o} } /
\sim_{ \hat{ \mathfrak{o} }^* \cap U } ) = Y $. To see this, let us prove
$ Y \subset \varphi^* ( \hat{ \mathfrak{o} } / \sim_{ \hat{ \mathfrak{o}
}^* \cap U } )$, since the other inclusion is certainly valid.
Take any $ ( (a_{\nu})^{\bullet}, (b_{\nu}) )^{\bullet} \in Y $, as $ \varphi^* $ is
surjective we can find $ (y_{\nu}) \in \mathbb{A}_f $ such that $ (
(a_{\nu})^{\bullet}, (b_{\nu}) )^{\bullet} = ( (1)^{\bullet}, (y_{\nu}) )^{\bullet} $.
Therefore, we can also find $ (z_{\nu}) $ in $ \hat{ \mathfrak{o} }^* $
such that $ ( (z_{\nu})^{\bullet}, (z_{\nu})^{-1} (y_{\nu}) ) = (
(a_{\nu})^{\bullet}, (b_{\nu}) ) $, and thus, $ (y_{\nu}) = (z_{\nu})
(b_{\nu}) \in \hat{ \mathfrak{o} } $ which means that $ ( (a_{\nu})^{\bullet},
(b_{\nu}) )^{\bullet} \in \varphi^* ( \hat{ \mathfrak{o} } / \sim_{ \hat{
\mathfrak{o} }^* \cap U } ) $.

This completes the proof for the case of no real places. If $K$ has one
real place, the proof will be exactly the same. But one should note that
in the computation above, one really needs the restriction to $ K^{\times}_{>0}
/ \mathfrak{o}^*_{>0} $ because for $ \lambda \in K $, $ (\lambda^{-1})^{\bullet}
\in \mathbb{A}_f^* $ lies in $ U $ if and only if $ \lambda $ is
positive.
\end{proof}

To get the relationship with our algebras $ \mathfrak{A} $, we remark that
there are canonical embeddings
\begin{equation}
C_0 ( \mathbb{A}_f / \sim_{ \hat{ \mathfrak{o} }^* \cap U } ) \rtimes
K^{\times} / \mathfrak{o}^* 
\hookrightarrow 
C_0 ( \mathbb{A}_f ) \rtimes K^\times \hookrightarrow C_0 ( \mathbb{A}_f ) \rtimes P_K 
\cong \mathfrak{A} (\mathfrak{o}) \nonumber
\end{equation}
if $K$ does not have real places and
\begin{equation}
C_0 ( \mathbb{A}_f / \sim_{ \hat{ \mathfrak{o} }^* \cap U } ) \rtimes
K^{\times}_{>0} / \mathfrak{o}^*_{>0} 
\hookrightarrow 
C_0 ( \mathbb{A}_f ) \rtimes K^{\times} \hookrightarrow C_0 ( \mathbb{A}_f ) \rtimes P_K 
\cong \mathfrak{A} (\mathfrak{o}) \nonumber
\end{equation}
for the case of one real place (see Theorem \ref{theo2} for the description of $ \mathfrak{A} $($ \mathfrak{o} $)).
Restricting these embeddings to the generalized Bost-Connes algebra $ \mathcal{A} $ (using the *-isomorphisms of Proposition \ref{algebraiso}), we get embeddings
$ \mathcal{A} \hookrightarrow \mathfrak{A}[\mathfrak{o}] $. At this point, we should note that there is no distinction between reduced or full crossed products as all the groups are amenable. Therefore, we really get embeddings.

\begin{remark}
For the case of purely imaginary number fields of class number one, there has been a remark in \cite{LLN} pointing in this direction (compare also the paper \cite{LaFr}).
\end{remark}

\subsection{Representation of the Bost-Connes Algebras}

From the point of view developed above, regarding the generalized Bost-Connes systems as subalgebras of our
algebras, it is possible to get an alternative description of $
\mathcal{A} $ as a universal $C^*$-algebra with generators and relations:

\begin{theorem}
\label{rep_BoCo}
Let $h_K = 1$ and assume that K has no real places. In
this case, $ \mathcal{A} $ is the universal $C^*$-algebra generated by 
nontrivial projections 
\begin{equation}
f(m,n) \text{ for every } m \in \mathfrak{o}^{\times} / \sim_{
\mathfrak{o}^* } \text{, } n \in \left( \mathfrak{o} / (m) \right)
/ \sim_{ \mathfrak{o}^* } \nonumber
\end{equation}
and isometries
\begin{equation}
s_p \text{ for each } p \in \mathfrak{o}^{\times} / \sim_{
\mathfrak{o}^* } \nonumber
\end{equation}
satisfying the relations
\begin{eqnarray}
&& s_p s_q = s_{pq} \text{ } \forall \text{ } p, q \in \mathfrak{o}^{\times} / \sim_{ \mathfrak{o}^* } \nonumber \\
&& f(1,0)=1 \nonumber \\
&& s_p f(m,n) s_p^* = f(mp,np) \text{ } \forall \text{ } p, m \in \mathfrak{o}^{\times} / \mathfrak{o}^* \text{, } 
n \in \left( \mathfrak{o} / (m) \right) / \sim_{ \mathfrak{o}^* } \nonumber \\
&& \sum_j f(mp,j) = f(p,k) \text{ } \forall \text{ } m, p \in \mathfrak{o}^{\times} / \mathfrak{o}^* \text{, } 
k \in \left( \mathfrak{o} / (p) \right) / \sim_{ \mathfrak{o}^* } \nonumber
\end{eqnarray}
where the sum in the last relation is taken over 
$ \pi_{m,mp}^{-1} (k) $ with the canonical projection 
\begin{equation}
\pi_{m,mp} : 
\left( \mathfrak{o} / (mp) \right) / \sim_{ \mathfrak{o}^* }
\longrightarrow \left( \mathfrak{o} / (p) \right) / \sim_{
\mathfrak{o}^* }. \nonumber
\end{equation}

If there is one real place, one has to substitute $ \mathfrak{o}^{\times}
$ by $ \mathfrak{o}^{\times}_{>0} $ and $ \mathfrak{o}^* $ by $
\mathfrak{o}^*_{>0} $.
\end{theorem}

Before we start with the proof, it should be noted that one can think of the
projection $ f(m,n) $ as $ \sum u^l e_m u^{-l} $ where the sum is taken
over all classes $ l + (m) $ in $ \mathfrak{o} / (m) $ which are in the
same equivalence class as $ n $ with respect to $ \sim_{\mathfrak{o}^* } $. 
This is exactly the way how these
elements are embedded into $ \mathfrak{A}[\mathfrak{o}] $. Moreover, using the
characteristic relations in $ \mathfrak{A} [\mathfrak{o}] $, the relations above can be
checked in a straightforward manner.

\begin{proof}
Let us prove the theorem in the case of no real places, the other case
can be proven in an analogous way.

The first step is to establish a *-isomorphism of the commutative
$ C^*$-algebras $ C( \hat{ \mathfrak{o} } / \sim_{ \hat{ \mathfrak{o} }^*
\cap U } ) $ and $ C^* ( \left\{ f(m,n) \text{ : } m \in \mathfrak{o}^{\times} / \sim_{
\mathfrak{o}^* } \text{, } n \in \left( \mathfrak{o} / (m) \right)
/ \sim_{ \mathfrak{o}^* } \right\} ) $.

We already had the description $ \hat{ \mathfrak{o} } = \plim \left\{ \mathfrak{o}
/ (m) \text{ ; } p_{m,lm} \right\} $. It will
be convenient to describe $ \hat{ \mathfrak{o} }^* \cap U = \overline{
\mathfrak{o}^* } $ in a similar way using projective limits. 

We claim: 
$ \hat{ \mathfrak{o}^* } \cap U \cong \plim \left\{ \left( \mathfrak{o}^* + (m) \right) /
(m) \text{ ; } p_{m,lm} \right\} $
where the $ p_{m,lm} $ are denoted as above, because again, they are
the canonical projections 

$ \left( \mathfrak{o}^* + (lm) \right) / (lm) \rightarrow \left( \mathfrak{o}^* + (m) \right) / (m) $.

To prove the claim, consider the following continuous embeddings
\begin{equation}
\mathfrak{o}^* \hookrightarrow \plim \left\{ \left( \mathfrak{o}^* + (m) \right) / (m)
\text{ ; } p_{m,lm} \right\} \hookrightarrow \hat{ \mathfrak{o} }. \nonumber
\end{equation}
Their composition is exactly the diagonal embedding of $ \mathfrak{o}^* $
into $ \hat{ \mathfrak{o} } $. 

Moreover, $ \plim \left\{ \left( \mathfrak{o}^* + (m) \right) / (m)
\text{ ; } p_{m,lm} \right\} $ (identified with its image in $
\hat{ \mathfrak{o} } $) is compact and contains $ \mathfrak{o}^* $. As it
follows from the construction of this projective limit that $
\mathfrak{o}^* $ (embedded in the projective limit) is dense, we have
proven the claim (compare Lemma \ref{calculationintersection}).

Furthermore, we have
\begin{eqnarray}
&& \hat{ \mathfrak{o} } / \sim_{ \hat{ \mathfrak{o} }^* \cap U }
\nonumber \\
& \cong & \plim \left\{ \mathfrak{o} / (m) \text{ ; } p_{m,lm} \right\}
/ \sim_{ \plim \left\{ \left( \mathfrak{o}^* + (m) \right) / (m) \text{ ; } p_{m,lm}
\right\} } \nonumber \\
& \cong & \plim \left\{ \left( \mathfrak{o} / (m) \right) \sim_{ \mathfrak{o}^* }
\text{ ; } p_{m,lm} \right\}. \nonumber
\end{eqnarray}
The first identification has already been proven; for the second one,
consider the following maps:
\begin{eqnarray}
\plim \left\{ \mathfrak{o} / (m) \right\} /
\sim_{ \plim \left\{ \left( \mathfrak{o}^* + (m) \right) / (m) \right\} } 
& \leftrightharpoons & 
\plim \left\{ \left( \mathfrak{o} / (m) \right) /
\sim_{ \mathfrak{o}^* } \right\} \nonumber \\
(a_m)^{\bullet} & \mapsto & (a_m^{\bullet}) \nonumber \\
(b_m)^{\bullet} & \mapsfrom & (b_m^{\bullet}) \nonumber
\end{eqnarray}

Existence and continuity of the upper map is given by the universal properties of
projective limits and quotient spaces. The lower map is well-defined
because of the following reason:

Let $ (b_m^{\bullet}) = (c_m^{\bullet}) $, we have to show $ (b_m)^{\bullet} = (c_m)^{\bullet}
$. 

For each $ m \in \mathfrak{o}^{\times} $, there is a $ r_m \in
\mathfrak{o}^* $ with the property that $ b_m + (m) = r_m c_m + (m) $.
But the net $ ( (r_m) )_m $ has a convergent subnet with limit $ (s_m)
$ in $ \plim \left\{ \left( \mathfrak{o}^* + (m) \right) / (m) \text{ ; } p_{m,lm}
\right\}  $ as this set is compact. And
the choice of the $ r_m $ ensures that we have 
$ (b_m) = (s_m) (c_m) $,
thus 
$ (b_m) \sim_{ \plim \left\{ \left( \mathfrak{o}^* + (m) \right) / (m) \text{ ; } p_{m,lm}
\right\} } (c_m) $. Therefore, the lower map exists as well.

Now, the upper map is a bijective continuous map, and the range is
Hausdorff, whereas the domain is quasi-compact. 
Hence it follows that these maps are mutually inverse
homoemorphisms.

After this step, we can now use Laca's result on crossed products by
semigroups to conclude the proof:

The universal $ C^* $-algebra with the generators and relations as
listed above is exactly given by the crossed product
\begin{equation}
C^* ( \left\{ f(m,n) \text{ : } m \in \mathfrak{o}^{\times} / \sim_{
\mathfrak{o}^* } \text{, } n \in \left( \mathfrak{o} / (m) \right)
/ \sim_{ \mathfrak{o}^* } \right\} ) \rtimes 
\left( \mathfrak{o}^{\times} / \sim_{ \mathfrak{o}^* } \right) \nonumber
\end{equation}
where we take the endomorphisms given by conjugation with $ s_p $. 
This is valid as both $ C^* $-algebras have the same universal properties. 

Furthermore, the identification above shows that 
\begin{eqnarray}
&& C^* ( \left\{ f(m,n) \text{ : } m \in \mathfrak{o}^{\times} / \sim_{
\mathfrak{o}^* } \text{, } n \in \left( \mathfrak{o} / (m) \right)
/ \sim_{ \mathfrak{o}^* } \right\} ) \rtimes \left( \mathfrak{o}^{\times} /
\sim_{ \mathfrak{o}^* } \right) \nonumber \\
& \cong &
C( \hat{ \mathfrak{o} } / \sim_{ \hat{ \mathfrak{o} }^*
\cap U } ) \rtimes \left( \mathfrak{o}^{\times} / \sim_{ \mathfrak{o}^* } \right) \nonumber
\end{eqnarray}
and the last $ C^* $-algebra is isomorphic to 
\begin{equation}
1_{ \hat{ \mathfrak{o} } / \sim_{ \hat{ \mathfrak{o} }^* \cap U } }
\left(
C_0 ( \mathbb{A}_f / \sim_{ \hat{ \mathfrak{o} }^* \cap U } ) \rtimes
K^{\times} / \mathfrak{o}^* 
\right)
1_{ \hat{ \mathfrak{o} } / \sim_{ \hat{ \mathfrak{o} }^* \cap U } } \nonumber
\end{equation}
by the work of Laca, since, in Laca's notation, 
$ 
C_0 ( \mathbb{A}_f / \sim_{ \hat{ \mathfrak{o} }^* \cap U } )
$
together with the action of
$
K^{\times} / \mathfrak{o}^*
$
is the minimal automorphism dilation corresponding to
$
C( \hat{ \mathfrak{o} } / \sim_{ \hat{ \mathfrak{o} }^*
\cap U } ) \rtimes \left( \mathfrak{o}^{\times} / \sim_{ \mathfrak{o}^* } \right) \nonumber
$
(see \cite{Lac}).

Finally, the result in Proposition \ref{algebraiso} gives us
\begin{eqnarray}
&& 1_{ \hat{ \mathfrak{o} } / \sim_{ \hat{ \mathfrak{o} }^* \cap U } }
\left(
C_0 ( \mathbb{A}_f / \sim_{ \hat{ \mathfrak{o} }^* \cap U } ) \rtimes
K^{\times} / \mathfrak{o}^* 
\right)
1_{ \hat{ \mathfrak{o} } / \sim_{ \hat{ \mathfrak{o} }^* \cap U } } \nonumber \\
& \cong & 1_Y \left( C_0 (X) \rtimes \mathbb{A}_f^* / \hat{ \mathfrak{o} }^* \right) 1_Y = \mathcal{A}. \nonumber
\end{eqnarray}

\end{proof}

\begin{remark}
Again, using $ \mathfrak{A}_r [ \mathfrak{o} ] \cong \mathfrak{A} [ \mathfrak{o} ] $, this result gives us a faithful (and again rather natural) representation of $ \mathcal{A} $ on $ \ell^2 ( \mathfrak{o} ) $.
\end{remark}

\begin{remark}
We can use Theorem \ref{rep_BoCo} to construct the extremal $KMS_\beta$-states of the C*-dynamical system 
$ \left( \mathcal{A}, \sigma_t \right) $, where $ \sigma_t (s_p) = N(p)^{it} s_p $ and $N$ is the norm in the number field $K$. Here, we view $ \mathcal{A} $ as the universal algebra as it is described in the previous Theorem. This is exactly the system considered in \cite{LLN}. 

We essentially follow the construction in \cite{BoCo}, THEOREM 25 (a), in the sense that for each element of the Galois group, we can construct a representation of $ \mathcal{A} $ using its universal property which yields the corresponding $KMS_\beta$-state.

First of all, we can associate to each $ \alpha \in Gal( K^{ab} / K ) $ the *-representation
\begin{equation}
\pi_\alpha : \mathcal{A} \longrightarrow \mathcal{L} \left( \ell^2 ( \mathfrak{o}^\times / \sim_{ \mathfrak{o}^* } ) \right) \nonumber
\end{equation}
by 
\begin{eqnarray}
&& \pi_\alpha (s_p) \xi_r = \xi_{pr} \nonumber \\
&& \pi_\alpha ( f(m,n) ) \xi_r 
= \begin{cases}
    \xi_r \text{ if } \overline{\alpha} ( r+(m) ) = n \in \left( \mathfrak{o}/(m) \right) / \sim_{ \mathfrak{o}^* } \\
    0 \text{ otherwise}.
  \end{cases} \nonumber
\end{eqnarray}
Here, $\overline{\alpha}$ is the image of $ \alpha $ under the composition
\begin{equation}
Gal(K^{ab}/K) \longrightarrow \hat{ \mathfrak{o} }^* / \overline{ \mathfrak{o}^* } \longrightarrow \left( \mathfrak{o} / (m) \right)^* / \sim_{ \mathfrak{o}^* }. \nonumber
\end{equation}

The existence of $\pi_{\alpha}$ follows from the universal property of $ \mathcal{A} $ described in Theorem \ref{rep_BoCo}.

Now, let us define $ H ( \xi_r ) = \log ( N(r) ) \xi_r $. With this operator we can construct the following $KMS_\beta$-state
\begin{equation}
\varphi_{\beta,\alpha} (x) = \zeta (\beta)^{-1} \text{tr} ( \pi_\alpha (x) e^{- \beta H} ) \nonumber
\end{equation}
where $\zeta$ is the zeta-function of the number field $K$.

This observation gives us candidates for the extremal $KMS_\beta$-states, and this construction follows an alternative, more operator-theoretic route compared to the rather measure-theoretic approach of \cite{LLN}. But it is another question whether these $ \varphi_{\beta, \alpha} $ are precisely all the extremal $KMS_\beta$-states for this C*-dynamical system, where $ 1 < \beta < \infty $. This is answered in the affirmative in \cite{LLN}, Theorem 2.1.

\end{remark}


\addcontentsline{toc}{section}{References}

\begin{thebibliography}{99}

\bibitem[Bla]{Bla} B. \textsc{Blackadar},
	\emph{Operator Algebras, Theory of C*-Algebras and von Neumann Algebras},
	Encyclopaedia of Mathematical Sciences, Vol. 122,
	Springer-Verlag, Berlin Heidelberg New York, 2006.

\bibitem[BoCo]{BoCo} J. B. \textsc{Bost} and A. \textsc{Connes},
  \emph{Hecke algebras, Type III Factors and Phase Transitions with Spontaneous Symmetry Breaking in Number Theory},
  Selecta Math., New Series, Vol. 1, \emph{3} (1995), 411-457.

\bibitem[CoMa]{CoMa} A. \textsc{Connes} and M. \textsc{Marcolli},
  \emph{From Physics to Number Theory via Noncommutative Geometry, PartI: Quantum Statistical Mechanics of 
  $\mathbb{Q}$-lattices} in Frontiers in Number Theory, Physics and Geometry, I, Springer-Verlag, Berlin Heidelberg New 
  York, 2006, 269-350.
  
\bibitem[CMR]{CMR} A. \textsc{Connes}, M. \textsc{Marcolli} and N. \textsc{Ramachandran},
  \emph{KMS states and complex multiplication},
  Selecta Math., New Series, \emph{11} (2005), 325-347.

\bibitem[Cun1]{Cun1} J. \textsc{Cuntz},
	\emph{C*-algebras associated with the $ax + b$-semigroup over $\mathbb{N}$}, 
	to appear in Proc. of the Conference on K-Theory and Non-Commutative Geometry, Valladolid 2006.

\bibitem[Cun2]{Cun2} J. \textsc{Cuntz},
	\emph{Simple C*-algebras generated by isometries},
	Comm. Math. Phys. \emph{57} (1977), 173-185. 

\bibitem[Lac]{Lac} M. \textsc{Laca},
  \emph{From endomorphisms to automorphisms and back: dilations and full corners},
  J. London Math. Soc. \emph{61} (2000), 893-904.

\bibitem[LaFr]{LaFr} M. \textsc{Laca} and M. \textsc{van Frankenhuijsen}
  \emph{Phase transitions on Hecke C*-algebras and class-field theory over $\mathbb{Q}$},
  J. reine angew. Math. \emph{595} (2006), 25-53.

\bibitem[LLN]{LLN} M. \textsc{Laca}, N. S. \textsc{Larsen} and S. \textsc{Neshveyev},
	\emph{On Bost-Connes type systems for number fields}, 
	preprint arXiv:0710.3452v2.

\bibitem[Neu]{Neu} J. \textsc{Neukirch},
	\emph{Algebraische Zahlentheorie},
	Springer-Verlag, Berlin Heidelberg New York, 1992.		

\bibitem[R\o r]{Ror} M. \textsc{R\o rdam},
	\emph{Classification of Nuclear C*-Algebras} in Classification of Nuclear C*-Algebras. Entropy in Operator Algebras,
	Encyclopaedia of Mathematical Sciences, Vol. 126,
	Springer-Verlag, Berlin Heidelberg New York, 2002.
	
\bibitem[Wei]{Wei} A. \textsc{Weil},
	\emph{Basic number theory} Die Grundlehren der mathematischen Wissenschaften,
	Band 144,	Springer-Verlag, Berlin Heidelberg New York, 1974.
	
\end{thebibliography}
\end{document}